\documentclass[a4paper,10pt]{article}
\usepackage[utf8]{inputenc}

\usepackage{amsfonts,latexsym,amstext}
\usepackage{amsmath}
\usepackage[active]{srcltx}
\usepackage{color}
\usepackage{amsmath,amsfonts,calrsfs,amssymb,color}
\usepackage{amssymb}
\newtheorem{theorem}{Theorem}[section]
\newtheorem{lemma}[theorem]{Lemma}
\newtheorem{proposition}[theorem]{Proposition}
\newtheorem{definition}{Definition}[section]

\newtheorem{hypothesis}[theorem]{Hypothesis}

\newtheorem{assumption}[theorem]{Assumption}
\newtheorem{remark}[theorem]{Remark}

\numberwithin{equation}{section}

\def\qed{{\hfill\hbox{\enspace${ \square}$}} \smallskip}
\def\sqr#1#2{{\vcenter{\vbox{\hrule height .#2pt \hbox{\vrule
 width .#2pt height#1pt \kern#1pt \vrule
width .#2pt} \hrule height .#2pt}}}}
\def\square{\mathchoice\sqr54\sqr54\sqr{4.1}3\sqr{3.5}3}

\def\ds{\begin{displaystyle}}
\def\eds{\end{displaystyle}}
\def\dis{\displaystyle }
\def\<{\langle }
\def\>{\rangle }

\def\dim{\noindent \hbox{{\bf Proof.} }}

\def\R{\mathbb R}
\def\N{\mathbb N}

\def\E{\mathbb E}
\def\P{\mathbb P}

\def\calb{{\cal B}}
\def\calc{{\cal C}}
\def\cald{{\cal D}}

\def\calf{{\cal F}}

\def\caln{{\cal N}}
\def\calp{{\cal P}}

\def\call{{\cal L}}

\title{}

\begin{document}

 \title{Partial smoothing of delay transition semigroups acting on special functions}
 \date{}
  \author{
 Federica Masiero, Gianmario Tessitore\\
 Dipartimento di Matematica e Applicazioni\\ Universit\`a di Milano Bicocca\\
 via Cozzi 53, 20125 Milano, Italy\\
 e-mail: federica.masiero@unimib.it, gianmario.tessitore@unimib.it}

\maketitle
\begin{abstract}
It is well known that the transition semigroup of an Ornstein Uhlenbeck process with delay is not strong Feller for small times, so it has no regularizing effects when acting on bounded and continuous functions. In this paper we study regularizing properties of this transition semigroup when acting on special functions of the past trajectory. With this regularizing property, we are able to prove existence and uniqueness of a mild solution for a special class of semilinear Kolmogorov equations; we apply these results to a stochastic optimal control problem.
\end{abstract}

%
%
 \section{Introduction}

In this paper we are at first concerned with the regularizing proerties of the transition semigroup related to the linear stochastic delay differential equation in $\R^n$
\begin{equation}
\left\{
\begin{array}
[c]{l}%
dy(t)  =a_0y(t)dt+\displaystyle \int_{-d}^0a_1(\theta)y(t+\theta)d\theta dt +\sigma dW_t
,\text{ \ \ \ }t\in[0,T] \\
y(0)=x_0\\
y(\theta)=x_1(\theta), \quad \theta \in [-d,0),
\end{array}
\right.  \label{eq-ritS_ritlin-intro}
\end{equation}
Due to the dependence at time $t$ on the past trajectory $(y_{t+\theta})_{\theta\in[-d,0]}$, the problem is intrinsicly infinite dimensional: it is reformulated in a space where both the evolution of the present and of the past trajectory are taken into account. For the sake of simplicity, in the introduction we mainly refer to the classical reformulation in the Hilbert space of square integrable past trajectories, see e.g. \cite{CM}, \cite{DM} and \cite{W}. Nevertheless in the paper we are able to consider the space $\cald$ where the past trajectoriy is a bounded cadlag function that has finite left limit, see Section \ref{sec-delayS} for more details. This will allow us to reduce the assumptions on the coefficients.
\newline The first Hilbert space in which we choose to reformulate  equation \eqref{eq-ritS_ritlin-intro}  is $H=\mathbb{R}^n\oplus L^{2}\left(\left[
-d,0\right] ,\mathbb{R}^n\right)$. For $t>0$  we introduce the operators
\[
e^{t A}:H\longrightarrow H,\quad e^{t A}\left(
\begin{array}{l}
x_{0}\\
x_{1}%
\end{array}
\right) =\left(
\begin{array}{l}
y\left( t\right) \\
y_{t}
\end{array}
\right), \text{ where } y_{t}\left(
\theta\right) =y\left( \theta+t\right),
\]
that define a $C_{0}$-semigroup with infinitesimal generator $A$, see (\ref{A}) for its definition. 
\newline Equation (\ref{eq-ritS_ritlin-intro}) can then be rewritten in an abstract way as a stochastic evolution equation in $H$:
\begin{equation}
\left\{
\begin{array}[c]{@{}l}%
dX(t)=AX(t)dt+GdW_{t},\text{ \ \ \ }%
t\in\left[ 0,T\right],\\
X_{0}=x,
\end{array}
\right. \label{eq-astrS-intro}%
\end{equation}
where
\[
G:\mathbb{R}^{n}\longrightarrow H,\text{ }G=\left(
\begin{array}[c]{@{}c@{}}%
\sigma\\
0
\end{array}
\right) ,\text{ }X_{t}=\left(
\begin{array}[c]{@{}c@{}}%
y\left( t\right) \\
y_{t}%
\end{array}
\right),\quad \text{and}\quad X_{0}=\left(
\begin{array}[c]{@{}c@{}}%
x_{0}\\
x_{1}%
\end{array}
\right),
\]
It is well known, see e.g. \cite{DP1} and \cite{DP2}, that for $t<d$ the Ornstein Uhlenbeck transition semigroup related to (\ref{eq-astrS-intro}) is not strong Feller, since in the associated deterministic controlled system
\begin{equation}
\left\{
\begin{array}[c]{@{}l}%
dz(t)=Az(t)dt+Gu(t)dt,\text{ \ \ \ }%
t\in\left[ 0,T\right],\\
z(0)=x,
\end{array}
\right. \label{sist-contr-det}%
\end{equation}
when $t<d$, it is not possible to steer to $0$ in time $t$ an initial state $x\neq 0$, even if $x\in \operatorname{Im}G$, see the discussion in \cite{Mas0}. So the transition semigroup related does not have any smoothing properties when acting on bounded and continuous function, even if we aim at achieving differentiability only in the present direction. 

The purpose of the present paper is to study the regularizing properties of the transition semigroup when acting on functions that have a special dependence on the past trajectoriy. The class of functionals considered here has already been studied in the literature, see e.g. \cite{FZ} and \cite{FRZ}; in particular we have been motivated by the functionals treated in \cite{RS}.
\newline Moreover in \cite{GoMa1}, for models arising in stochastic optimal control problems with delay in the control,  the regularizing properties of the transition semigroup are studied in suitable directions and on a class of special functions that arise naturally in that case.

Coming into more details, the class of functionals we are concerned with here is constructed as follows: consider a Borel measurable and bounded function $\bar\phi:\R^n\rightarrow\R$ and define  $\phi=\bar\phi\circ P$ by:
\begin{equation}\label{fiS-intro}
\phi(x)=\bar\phi(P(x)) \quad
\forall\, x=\left( \begin{array}{l}
x_0 \\x_1
\end{array}
 \right)\in H,\quad P:H\rightarrow \R^n, \quad P\left(\begin{array}{l}
                             x_0\\
                             x_1\end{array}\right)= \alpha_0x_0 +\int_{-d}^0f(\theta)x_1(\theta)d\theta,
\end{equation}
where $\alpha_0\in\operatorname{Mat}(n\times n)$ and $f\in L^2([-d,0],\operatorname{Mat}(n\times n))$.
\newline Under suitable assumptions on $\alpha_0$ and $f$, we are able to prove that, if $\bar\phi$ is a bounded continuous function and $\phi$ is defined as in (\ref{fiS-intro}), then the Ornstein Uhlencbeck transition semigroup $R_t,\,t>0, $ maps $\phi $ into a differentiable function
and 
\begin{equation}\label{stima_der-intro}
  \vert\nabla R_t[\phi](x)h\vert \leq C \frac{\vert h\vert}{\sqrt{t}}\Vert \phi\Vert_\infty.
 \end{equation}
Then, see Section 3, we consider the process $X$  living in the smaller space $\cald$ where the past trajectory is continuous apart from a finite number of points; this enables us to consider more general functions 
\begin{equation}\label{fi-gen-intro}
\phi(x)=\bar\phi(\calp(x)) \quad
\forall x=\left( \begin{array}{l}
x_0 \\x_1
\end{array}
 \right)\in \cald,\quad
 \calp\left(\begin{array}{l}
                             x_0\\
                             x_1
\end{array}\right)=
                  \mu(\left\lbrace0\right\rbrace)x_0+   \int_{[-d,0)}\mu(d\theta)x_1(\theta)  ,
\end{equation}
where $\mu $ is a finite regular measure on $[-d,0]$, see Hypothesis \ref{ip-mu-stima} for more details.
Moreover we are able to prove that the regularizing properties hold true also for suitable perturbationns of the Ornstein Uhlenbeck process (\ref{eq-astrS-intro}).

Such regularizing properties allow us to solve, by generalizing a fixed point argument, semilinear Kolmogorov equations as
\begin{equation}\label{Kolmo-formale0-intro}
  \left\{\begin{array}{l}\dis
-\frac{\partial v(t,x)}{\partial t}=\call[v(t,\cdot)](x) +
\psi (v(t,x), \nabla v(t,x)G),\qquad t\in [0,T],\,
x\in \cald,\\
\\
\dis v(T,x)=\phi(x),
\end{array}\right.
\end{equation}
where $\phi$ is a bounded and continuous function with the structure described in (\ref{fiS-intro}) ( or in (\ref{fi-gen-intro})), 
and $\call$ is the generator of the Ornstein-Uhlenbeck semigroup. We only mention that, in Section \ref{subsec-Kolmo}, we preliminarly study linear Kolmogorov equations, that is equations \eqref{Kolmo-formale0-intro} with $\psi \equiv 0$.
\newline By applying formally the variation of constants formula, see e.g. \cite{DP3}, the semilinear Kolmogorov equation (\ref{Kolmo-formale0}) can then be rewritten in its ``mild formulation''
\begin{equation}
v(t,x) =R_{T-t}[\phi](x)+\int_t^T R_{s-t}\left[
\psi(v(s,\cdot),\nabla v(s,\cdot)G)\right](x)\; ds,\qquad t\in [0,T],\
x\in \cald,\label{solmildKolmo_sem-intro}
\end{equation}
Semilinear Kolmogorov equations related to problem with delay in the state are solved in \cite{FMT} by using a probabilistic approach based on backward stochastic differential equations under differentiability assumptions on the coefficients of the equations, in particular on $\phi$. Here, due to the regularizing properties of the transition semigroup, we are able to require only continuity of the final datum. Moreover when the coefficients are assumed to be twice differentiable, semilinear Kolmogorv equations like \eqref{solmildKolmo_sem-intro} are solved, in classical sense, in \cite{MOTZ}, generalizing the linear case treated in \cite{FZ}. Besides these papers on mild and classical solutions, for which differentiability of the solution is required, path dependet PDEs are largely studied in the literature, mainly in the viscosity sense. In this framework a more general structure than \eqref{Kolmo-formale0-intro} can be considered, see e.g. \cite{EKTZ} and  \cite{RR}, see also \cite{CGRST} where infinite dimensional path dependent PDEs are considered.

In the present paper, the results on the existence of a mild solution of equation (\ref{Kolmo-formale0-intro}) are applied to a stochastic optimal control problem where to the controlled  equation 
\begin{equation}
\label{eq:SDEx:control-intro}
\left\{
\begin{array}
[c]{l}%
dy(t)  =a_0y(t)dt+\displaystyle \int_{-d}^0a_1(\theta)y(t+\theta)d\theta dt +\sigma u(t)dt+\sigma dW_t
,\text{ \ \ \ }t\in[0,T] \\
y(0)=x_0\\
y(\theta)=x_1(\theta), \quad \theta \in [-d,0),
\end{array}
\right. 
\end{equation}
we associate the cost functional $J$:
\begin{equation}
J\left(  t,x,u\right)  =\mathbb{E}\int_{t}^{T} g\left(u(s)\right)ds+\mathbb{E}\bar\phi\left(\calp y^{u}_T\right). 
\label{cost-intro}%
\end{equation}
Notice the dependence on the whole trajectory in the final cost, where $\calp $ is the kind of dependence on the past that we are able to handle reformulating the problem in $\cald$. On the other side notice that the running cost does not depend on $y$.
\newline The related Hamilton Jacobi Bellamn equation has the structure of equation (\ref{Kolmo-formale0-intro}) and it turns out that its unique mild solution is the value function of the control problem. Moreover the optimal control  can be characterized in a feeedback form, that involves the gradient of the mild solution $v$ to the HJB equation. We notice that first order regularity of $v$ is guaranteed by our approach. 


The paper is organized as follows in Sections \ref{sec-delayS} we study the regularizing properties of the transition semigroup respectively in the Hilbert space $H$ and in the Banach space $\cald$, while  in \ref{sec:reg_OUpert}  we show how these regularizing properties are inherited by the perturbed Ornstein Uhlenbeck processes,  both in the Hilbert space $H$ and in the Banach space $\cald$. In Section \ref{sec-mildKolmo} the results are applied to the solution in mild sense of linear and semilinear Komogorv equations, finally in Section \ref{sec-control} applications to control are given.

\section{Partial smoothing of Ornstein-Uhkenbeck transition semigroups}
\label{sec-delayS}

In a complete probability space $(\Omega, \calf,  \P)$ we consider the following controlled stochastic differential equation in $\R^n$ with delay in the state:
\begin{equation}
\left\{
\begin{array}
[c]{l}%
dy(t)  =a_0y(t)dt+\displaystyle \int_{-d}^0y(t+\theta)a_1(d\theta) dt 
+\sigma dW_t
,\text{ \ \ \ }t\in[0,T] \\
y(0)=x_0\\
y(\theta)=x_1(\theta), \quad \theta \in [-d,0)\,\,a.e.,
\end{array}
\right.  \label{eq-ritS_ritlin}
\end{equation}
%
where $a_0\in \operatorname{Mat}(n\times n)$ and $a_1$ is an $n\times n$ matrix valued finite regular measure, and it is such that 
\begin{equation}\label{ip-eta-stima}
a_1({0})=0,
\end{equation}
where here by $0$ we mean the $n\times n$ matrix identically equal to $0$.
We notice that \eqref{ip-eta-stima} implies that
$$
\int_{-t}^0 a_1(d\theta)\rightarrow 0 \quad \text{ as }t\rightarrow 0.
$$
The process $ W$ in \eqref{eq-ritS_ritlin} is a standard Wiener process in $\mathbb{R}^n$ and $\sigma \in \operatorname{Mat}(n\times n)$. The
value $d>0$ denotes the maximum delay, $x_0\in\R^n$, $x_{1}\in L^{2}([-d,0],\mathbb{R}^n)$; moreover in the following $y_t$ denotes the past trajectory from time $t-d$ to time $t$:
$$
y_t(\theta):=y(t+\theta), \;\theta \in [-d,0],\quad y_t\in L^{2}([-d,0]
,\mathbb{R}^n).
$$
%
\newline
Define the Hilbert space $H=\mathbb{R}^n\oplus L^{2}\left(\left[-d,0\right] ,\mathbb{R}^n\right)$. 
Let $y(t)$ be the solution at time $t>0$ of the following linear deterministic problem with delay
\begin{equation}\label{ed-delay-det}
\left\{
\begin{array}
[c]{l}%
dy(t)  =a_0y(t)dt+\displaystyle \int_{-d}^0y(t+\theta)a_1(d\theta)dt
,\text{ \ \ \ }t\in[0,T] \\
y(0)=x_0\\
y(\theta)=x_1(\theta), \quad \theta \in [-d,0).,\,a.e.
\end{array}
\right.  
\end{equation}
It turns out that
\begin{equation}\label{semigr-delay}
e^{t A}:H\longrightarrow H,\quad e^{t A}\left(
\begin{array}[c]{@{}c@{}}%
x_{0}\\
x_{1}%
\end{array}
\right) =\left(
\begin{array}[c]{@{}c@{}}%
y\left( t\right) \\
y_{t}%
\end{array}
\right),\quad t\geq 0 \text{ with } y_{t}\left(
\theta\right) =y\left( \theta+t\right)
\end{equation}
defines a $C_{0}$-semigroup in $H$; see, e.g.,~\cite{DM} and \cite{W}.
The infinitesimal generator $A$ of $(e^{t A})_{t\geq 0}$, is
given by
\begin{equation}\label{A}
\begin{array}[c]{l}%
\mathcal{D}\left( A\right) =\left\{ \left(
\begin{array}[c]{@{}c@{}}%
x_{0}\\

x_{1}%
\end{array}
\right) \in H,x_{1}\in H^{1}\left( \left[ -d,0\right] ,\mathbb{R}%
^{n}\right) ,x_{1}\left( 0\right) =x_{0}\right\},\\
\\

Ah=A\left(
\begin{array}[c]{@{}c@{}}%
x_{0}\\

x_{1}%
\end{array}
\right) =\left(
\begin{array}[c]{@{}c@{}}%
a_0 x_0+
{\displaystyle\int_{-d}^{0}}
x_{1}\left( \theta\right) a_1\left( d\theta\right) \\

dx_{1}/d\theta
\end{array}
\right) .
\end{array}
\end{equation}
By setting%
\begin{equation}\label{G}
G:\mathbb{R}^{n}\longrightarrow H,\text{ }G=\left(
\begin{array}[c]{@{}c@{}}%
\sigma\\
0
\end{array}
\right) ,\text{ }X_{t}=\left(
\begin{array}[c]{@{}c@{}}%
y\left( t\right) \\
y_{t}%
\end{array}
\right),\quad \text{and}\quad X_{0}=\left(
\begin{array}[c]{@{}c@{}}%
x_{0}\\
x_{1}%
\end{array}
\right),
\end{equation}
problem (\ref{eq-ritS_ritlin}) can be rewritten in an abstract way as an Ornstein-Uhlenbeck process in $H $ given by
\begin{equation}
\left\lbrace \begin{array}
[c]{l}
dZ^x(t)  =AZ^x(t) dt+GdW_t
,\text{ \ \ \ }t\in[ 0,T] \\
Z^x(0)  =x=\left(\begin{array}{l}x_0\\x_1\end{array}\right).
\end{array}
\right.   \label{eq-astrS-OU}%
\end{equation}
Taking the integral mild form of (\ref{eq-astrS-OU}) we have
\begin{equation}
Z^x(t)  =e^{tA}x+\int_0^te^{(t-s)A}GdW_s
,\text{ \ \ \ }t\in[ 0,T]. \\
  \label{eq-astr-mild-OU}%
\end{equation}
The Ornstein-Uhlenbeck transition semigroup $R_t$ is defined by setting,
for every measurable function $f:H\rightarrow\R$ and for every $x\in H$,
\begin{equation}
 \label{ornstein-semS}
R_t[f](x)=\E [f(Z^x(t))]=\int_H f(z)\caln(e^{tA}x,Q_t)(dz)
=\int_H f(z+e^{tA}x)\caln(0,Q_t)(dz).
\end{equation}
We look for partial smoothing properties of  the transition semigroup $(R_t)_{t>0}$ when acting on special functions defined as follows.
We consider the map
\begin{equation}\label{P}
 P:H\rightarrow \R^n, \quad P\left(\begin{array}{l}
                             x_0\\
                             x_1\end{array}\right)= \alpha_0x_0 +\int_{-d}^0f(\theta)x_1(\theta)d\theta
\end{equation}
where $\alpha_0\in\operatorname{Mat}(n\times n)$ and $f\in L^2([-d,0],\operatorname{Mat}(n\times n))$.
\newline Given any Borel measurable and bounded function
$\bar\phi:\R^n\rightarrow\R$ we define $\phi :H\rightarrow\R$ by setting
\begin{equation}\label{fiS}
\phi(x)=\bar\phi(P(x)) \quad
\forall x=\left( \begin{array}{l}
x_0 \\x_1
\end{array}
 \right)\in H,
 \end{equation}
so that $\phi=\bar\phi\circ P$.

In the following we prove that the transition semigroup maps bounded and measurable function defined according to \eqref{fiS} into differentiable ones.
We also underline the fact that it is well known, see e.g. \cite{DP2}, that if $t>d$, the transition semigroup $R_t$ is strong Feller: it maps bounded and measurable functions into differentiable ones. Here we are able to extend this regularizing property also for $t<d$, but only when the transition semigroup $R_t$ acts on special functions like the ones defined in (\ref{fiS}).
This regularizing property is  the basic tool to prove existence of a solution to the Kolmogorov equations that we study in Section \ref{sec-mildKolmo}.

In order to prove this partial smoothing property of the transition semigroup, we have to make some assumptions, included the invertibility of the diffusion coefficient $\sigma$ in equation \eqref{eq-ritS_ritlin}.
\begin{assumption}\label{ass} We will work under one of the following non-degeneracy assumptions  (mainly the first)
\begin{equation} \label{ass-1}\tag{A1} \hbox{ both $n\times n$ matrices $\sigma$ and $\alpha_0 $ are invertible.}
\end{equation}
\begin{equation}  \label{ass-2}\tag{A2} \hbox{ $\sigma$ is invertible, $\alpha_0=0$, there exists  an invertible matrix $f_0$ such that $s^{-1}\dis\int_{-s}^0 f(\theta) d \theta \rightarrow f_0 $  in $\R^{n\times n }$}
\end{equation}
\end{assumption}
\begin{theorem}\label{prop-reg}
Let $(R_t)_{t\geq 0}$ be the transition semigroup  related  to equation (\ref{eq-astrS-OU}) and defined accordingly to formula (\ref{ornstein-semS}). 
Let $\phi$ be a bounded and measurable function defined as in (\ref{fiS}) with 
$\bar\phi \in B_b(\R^n)$.
Then $R_t[\phi]:H\rightarrow\R$ is differentiable. Moreover, the gradient $\nabla R_t[\phi]$ can be estimated depending on the non degeneracy assumption we choose. Namely there exists $C>0$ such that for
 all   $h\in H$
 \begin{equation}\label{stima_der}
 \hbox{ if \eqref{ass-1} holds then } \vert\nabla R_t[\phi](x)h\vert \leq C \frac{\vert h\vert}{\sqrt{t}}\Vert \phi\Vert_\infty;
 \end{equation}
  \begin{equation}\label{stima_der_bis}
 \hbox{ if \eqref{ass-2} holds then } \vert\nabla R_t[\phi](x)h\vert \leq C \frac{\vert h\vert}{t}\Vert \phi\Vert_\infty.
 \end{equation}
\end{theorem}
\dim We start by proving the strong Feller property and the related estimate (\ref{stima_der}). We compute, for any $h\in H$
\begin{align}\label{Rt-lip}
& R_t[\phi](x+h)- R_t[\phi](x)\\  \nonumber
&=\int_{H}\bar\phi(P(y+e^{tA}(x+h)))\caln(0,Q_t)(dy)-\int_{\R}\bar\phi(P(y+e^{tA}x))\caln(0,Q_t)(dy)\\ \nonumber
&=\int_{\R^n}\bar\phi(z+Pe^{tA}(x+h))\caln(0,PQ_tP^*)(dz)-\int_{\R^n}\bar\phi(z+Pe^{tA}x)\caln(0,PQ_tP^*)(dz)\\ \nonumber
&=\int_{\R^n}\bar\phi(z+Pe^{tA}x)\caln(Pe^{tA}h,PQ_tP^*)(dz)-\int_{\R^n}\bar\phi(z+Pe^{tA}x)\caln(0,PQ_tP^*)(dz),
\end{align}
We have to estimate from below the covariance operator $PQ_tP^*$: for $\xi\in \mathbb{R}^n$
\begin{equation}\label{eq:cov}\<\xi,PQ_tP^*\xi\>_{\R^n}=\int_0^t \<\xi,Pe^{sA}GG^*e^{sA^*}P^*\xi\>_{\R^n}\,ds=\int_0^t |G^*e^{sA^*}P^*\xi|_{\R^n}^2	,ds=\int_0^t \left[\sup_{ |\eta|=1}\< Pe^{sA}G\eta,\xi\>_{\R^n}^2\right]ds
\end{equation}
Fixed $\eta \in \mathbb{R}^n$ with $|\eta|=1$ we have, for $t<d$:
\begin{equation}\label{eq:P} Pe^{t A}G\eta=\alpha_0 \bar y(t)+\int_{-t}^0 f(\theta)\bar y(t+\theta)d\theta
\end{equation}
where $\begin{pmatrix}
\bar y(t)\\
\bar y_t
\end{pmatrix} = e^{tA} \begin{pmatrix}
\sigma \eta\\
0
\end{pmatrix}$. By a straight-forward application of the variation of constants formula, again for $t<d$
$$\bar y(t)=e^{t a_0}(\sigma \eta)+\int_0^t e^{(t-s)a_0} \int_{-s}^0 \bar y(s+\theta)a_1(d\theta) \,ds=e^{t a_0}(\sigma \eta)+\int_{-t}^0\int_{-\theta}^t e^{(t-s)a_0}  \bar y(s+\theta)\,ds \,a_1(d\theta)  $$
By standard estimates, also recalling that $\bar y_s(\cdot) $ is square integrable and $a_1$ is a finite measure (with finite total variation) we get
\begin{equation}
\label{stima-bar y}
\bar{y}(t)-e^{t a_0}(\sigma \eta):=r_0(t)\; \hbox{ with }\;|r_0(t)|\leq c |a_1|([-t,0])\vert y\vert_{L^2([0,t])}t^{1/2},
\end{equation}
where by $|a_1|([-t,0])$ we denote the total variation of $a_1$ on the interval $[-t,0]$.
Noting that $\vert y\vert_{L^2([0,t])}\rightarrow 0$ as $t\rightarrow 0 $, we can deduce that 
\begin{equation}
\label{stima-bar y-o}
t^{-1/2}\sup_{s\in[0,t]}\left\vert\bar{y}(s)-e^{s a_0}(\sigma \eta)\right\vert\text{ as }t\searrow 0.
\end{equation}
Thus
$$ Pe^{t A}G\eta=\alpha_0 [e^{t a_0}(\sigma \eta)+r_0(t)]+\int_{-t}^0 f(\theta)[e^{(t+\theta) a_0}(\sigma \eta)+r_0(t+\theta)]d\theta=\alpha_0 e^{t a_0}(\sigma \eta)+\int_{-t}^0 f(\theta)(\sigma \eta)d\theta+r_1(t)$$
where under \eqref{ass-1} $t^{-1/2}r_1(t) \rightarrow 0 $ as $t\searrow 0$ and under \eqref{ass-2} $t^{-3/2}r_1(t) \rightarrow 0 $ as $t\searrow 0$.
Thus under \eqref{ass-1}:
$$Pe^{t A}G\eta=\alpha_0 \sigma \eta +r_2(t) \, \hbox{ where  }\; t^{-1/2}r_2(t)\rightarrow 0\text{ as }t\searrow 0,$$
and under \eqref{ass-2}:
$$Pe^{t A}G\eta=t f_0 \sigma \eta +r_3(t) \, \hbox{ where  }\; t^{-1}r_3(t)\rightarrow 0 \text{ as }t\searrow 0.$$
Taking into account \eqref{eq:cov}, \eqref{eq:P} invertibility of $\alpha_0 \sigma $ and $f_0 \sigma$ we have, for $t\in (0,\epsilon)$ and $\epsilon$ small enough :
$$\<\xi,PQ_tP^*\xi\>_{\R^n}\geq c |\xi|^2 t^{\gamma}$$
with $\gamma=1$ under \eqref{ass-1} and $\gamma=3$ under \eqref{ass-2}.
To shorten notation in the following we set  $
\bar Q_t: =PQ_tP^*.
$

\medskip

Coming back to (\ref{Rt-lip}), it turns out that for any $\phi\in B_b(\R^n)$, $R_t[\phi]$ is continuous since it is Lipschitz continuous moreover for any $h \in H$, $|h|=1$, setting $y=Pe^{tA}h$ we get for $t>0
$:
 \begin{align*}
 &\< \nabla R_t[\phi],h\>\\
 & =\lim_{s\rightarrow 0}\dfrac{1}{s}\left(\int_{\R^n}\bar\phi(z+Pe^{tA}x)\caln(s Pe^{tA}h,\bar Q_t)(dz)
  -\int_{\R^n}\bar\phi(z+Pe^{tA}x)\caln(0,\bar Q_t)(dz)\right)\\
   &=\lim_{s\rightarrow 0}\dfrac{1}{s}\int_{\R^n}\bar\phi\left(z+Pe^{tA}x\right)\left(
1-e^{ s \left\langle \bar Q_{t}^{-1/2}
Pe^{tA}h,\bar Q_{t}^{-1/2}z\right\rangle_{\R^n}
-\frac{s^2}{2}\left|\bar Q_{t}^{-1/2} Pe^{tA}h\right|_{\R^n}^{2}} \right)\caln(0,\bar Q_t)(dz)\\ &
= \int_{\R^n}\bar\phi(z+Pe^{tA}x)\left\langle\bar Q_{t}^{-1/2}
Pe^{tA}h,\bar Q_{t}^{-1/2}z\right\rangle_{\R^n}\caln(0,\bar Q_t)(dz)
\\ &
= \int_{\R^n}\bar\phi(\bar Q_{t}^{1/2}\zeta+Pe^{tA}x)\left\langle\bar Q_{t}^{-1/2}
Pe^{tA}h, \zeta\right\rangle_{\R^n}\caln(0, I)(d\zeta)
  \end{align*}
and we can conclude that:
\begin{equation*}
 \vert \< \nabla R_t[\phi],h\>\vert \leq {C}{t^{-\frac{\gamma}{2}\gamma}}\Vert \bar\phi \Vert_\infty = {C}{t^{-\frac{\gamma}{2}}}\Vert \phi \Vert_\infty
\end{equation*}
where $\gamma=1$ under \eqref{ass-1} and $\gamma=3$ under \eqref{ass-2}

\qed

\smallskip

When solving the HJB equation we will need also the following result, which turns out to be a generalization of the previous smoothing result. Indeed in the previous Theorem the proof is based on the invertibility of the operator $\bar Q_t$ and on the estimate of $\vert\bar Q_t\vert^{-1/2} $; in the following proposition we consider the operator
\begin{equation}\label{bar-Qst}
\bar Q^s_t:=Pe^{sA}Q_{t-s}e^{sA^*}P^*=\int_s^t Pe^{rA}GG^*e^{rA^*}P^*\,dr,\quad 0\leq s <t\leq T.
\end{equation}
We present a result only under \eqref{ass-1}, an analogous result under \eqref{ass-2} folllows in a similar way.
\begin{proposition}\label{prop_barQ}
 Let $A,\,G$ and $P$ be defined respectively in \eqref{A}, \eqref{G} and \eqref{P} and let \eqref{ass-1} holds true. Let $\bar Q^s_t$ be defined in \eqref{bar-Qst}, $0\leq s<t\leq d$. Then there exists $\bar t$ such that $\forall\, 0\leq s<t\leq \bar t$, $\forall\,\xi\in \R^n$
\begin{equation}\label{eq:stima_barQ_prel}
\<\xi,P\bar Q^s_tP^*\xi\>_{\R^n}\geq c|\xi|^2 (t-s) 
\end{equation}
so that 
\begin{equation}\label{eq:stima_barQ}
 \left(Q^s_t\right)^{-1/2 }\vert \leq c (t-s)^{-1/2}.
\end{equation}
\end{proposition}
\dim We evaluate, arguing in a similar way to what done in the proof of Theorem \ref{prop-reg}, 
$$
\<\xi,\bar Q^s_t\xi\>_{\R^n}=\int_s^t\vert G^*e^{rA^*}P^*\xi\vert_{\R^n}\,dr=\int_s^t\sup_{\vert \eta\vert =1}\<P e^{rA}G\eta,\xi\>_{\R^n}\,dr.
$$
We recall that, for $r<d$,
$$
P e^{rA}G\eta=\alpha_0\sigma \eta+r_2(r), \text{ where } r_2(r)=o(r^{1/2}),
$$
thus $\forall\, \delta >0$ there exists $t_\delta$ such that for $\xi\in\R^n$
\[
\<P e^{rA}G\eta,\xi\>\geq \<\alpha_0\sigma\eta,\xi\>-\delta|\xi|,\; \forall \, r\in [0,t_\delta].
\]
Choosing $\bar\eta= \dfrac{(\alpha_0\sigma)^{-1}\xi}{\vert (\alpha_0\sigma)^{-1}\xi\vert}$ we get 
$$
\<(\alpha_0\sigma)\bar\eta,\xi\>\geq \vert(\alpha_0\sigma)\vert \vert \xi\vert
$$
and finally choosing $\delta =\frac{1}{2}\vert (\alpha_0\sigma)\vert$ we get 
$$
\<P e^{rA}G\bar \eta,\xi\>\geq \frac{1}{2}\vert (\alpha_0\sigma)\vert\vert \xi\vert,\, \forall r\in[0.t_\delta],
$$ 
and the claim follows.
%
\qed

\begin{remark}\label{rm:no_riduzione}
We have focused our attention on regularizing properties of the transition semigroup when acting on functions 
defined as in $\eqref{fiS}$; we would like to stress the fact that also in this case the model doesn't allow a simpler reduction. Indeed letting $y$ be the solution of \eqref{eq-ritS_ritlin-intro} and setting 
\begin{equation}\label{Y}
 Y(t):=\alpha_0y(t)+\int_{-d}^0 y(t+\theta)f(\theta)\,d\theta.
\end{equation} one heuristically obtains (for regular $f$)
$$d_t Y(t)=\alpha_0 dy(t)+f(0)y(t)-f(-d)y(t-d)-\int_{-d}^0 f'(\theta)y(t+\theta)\,d\theta$$
It is then clear that (even when $a_1\equiv0$ and consequently we just have $dy(t)=a_0y(t) dt + \sigma dw_t$) our model does not give rise to a markovian dynamic neither for $Y$ nor for $(y,Y)$.

\end{remark}

We now consider the case when the past trajectory is a cadlag function. Indeed, 
if in equation \eqref{eq-ritS_ritlin} the initial past trajectory $x_1$ is a continuous functions, the pair $\Big( \begin{array}{l}
y(t)\\y_t
\end{array}\Big)$ evolves in a Banach space $\cald$, continuously and densely embedded in $H$, and that we are going to introduce.

We denote by $D_b([-d,0),\R^n)$ the set of bounded cadlag functions that have finite left limit for $\theta\nearrow 0$ and we define the product space 
\begin{equation}\label{cald}
\cald:=\left\lbrace x=\left(\begin{array}{l}x_0\\x_1
\end{array}
\right)\in \R^n\times D_b([-d,0), \R^n )
\right\rbrace.
\end{equation}
The space $\cald $ turns out to be a Banach space if it is endowed with the norm
\[
\left\| \left(\begin{array}{l}x_0\\x_1
\end{array}
\right)\right\|_\cald:=\vert x_0\vert+\Vert x_1\Vert_\infty.
\]
It turns out that if in equation (\ref{eq-ritS_ritlin}) $x=\left(\begin{array}{l}x_0\\x_1
\end{array}
\right)\in\cald$, see e.g. \cite{moha}, then for all $t>0$ $\left(\begin{array}{l}   
y(t)\\
y_t
\end{array}\right)\in\cald.
$
Moreover if for $t\geq 0$ we set $e^{t A}:H\longrightarrow H,\quad e^{t A}\left(
\begin{array}[c]{@{}c@{}}%
x_{0}\\
x_{1}%
\end{array}
\right) =\left(
\begin{array}[c]{@{}c@{}}%
y\left( t\right) \\
y_{t}%
\end{array}
\right), \text{ with } y_{t}\left(
\theta\right) =y\left( \theta+t\right)$ then $e^{tA}$ is the generator of a semigroup in $\cald$, which is the restriction to $\cald$ of the semigroup in $H$. Moreover, see \cite{FZ}, for some constant $C>0$
$$
\vert e^{tA}\vert_{\call(\cald,\cald)}\leq C,\; t\in[0,T].
$$ 
Following again \cite{FZ} , whenever $x\in\cald$, equation (\ref{eq-ritS_ritlin}) can be reformulated as an Ornstein Uhlenbeck process in $\cald$, in the sense that letting again $Z^x(t)=\Big( \begin{array}{l}
y(t)\\y_t
\end{array}\Big)$
\begin{equation}
Z^x(t)  =e^{tA}x+\int_0^te^{(t-s)A}GdW_s
,\text{ \ \ \ }t\in[ 0,T], \\
  \label{eq-astr-mild-OU-cald}%
\end{equation}
and now we remark that if $x\in \cald$ $e^{tA}x\in \cald$,
and the operator $G$ this time is the inclusion of $\R^n$ into the product space $\cald$:
\begin{equation}\label{G-gen}
 G:\R^n\rightarrow \cald,\; y\mapsto \left( \begin{array}{l}y\\0   
 \end{array}\right).
\end{equation}
We notice that being $\cald$ a Banach space lacking the topological properties needed to build in an infinite dimensional stochastic calculus, it is not even obvious how to define the stochastic convolution
\begin{equation*}
  W^r_A(t)=\int_r^te^{(t-s)A}G\,d W_s\ ,\quad 0\leq r\leq t\leq T,
\end{equation*}
(that we simply denote by $W_A(t)$ if $r=0$).
The construction of $W^r_A(t)$ together with its properties have been stated in \cite{FZ}
In particular in \cite{FZ} it has been proved that $W^r_A$ is gaussian and has continuous trajectories in the product space $\R^n\times E$, where $E=\left\lbrace f\in C([-d,0),\R^n) :\exists\, \lim_{r\nearrow 0} f(r)  \right\rbrace$. Clearly $\R^n\times E\subset \cald$, so the Ornstein-Uhlenbeck process $Z^x$ solution to equation \eqref{eq-astr-mild-OU-cald} is a well defined process with values in $\cald$.

The associated Ornstein-Uhlenbeck transition semigroup $R_t$, is defined by setting, for every measurable function $f:\cald\rightarrow\R$, and for every $x\in \cald$,
\begin{equation}
 \label{ornstein-semS-gen}
R_t[f](x)=\E [f(Z^x(t))]
\end{equation}

We aim at proving that the transition semiproup $(R_t)_{t>0}$ is regularizing in $\cald$ when acting on special functions similarly to the regularizing properties we have proved in $H$.
\newline When working in $\cald$, the class of special functions we can consider is larger. In order to introduce this class of special functions, we consider a second regular measure $\mu$ satisfying the following:
\begin{hypothesis}\label{ip-mu}
The measure $\mu$ is a regular measure on $[-d,0]$ with values in $\operatorname{Mat}(n\times n)$ and it is such that 
\begin{equation}\label{ip-mu-stima}
\alpha_0:=\mu(\left\lbrace0\right\rbrace)>0.
\end{equation}
 \end{hypothesis}
It follows that for all $A\subset[-d,0]$, $\mu(A)=\alpha_{0}\delta_0(A)+\bar\mu(A\setminus\left\lbrace 0\right\rbrace)$ where $\bar \mu $ is a regular measure on $[-d,0]$ with $\bar\mu (\lbrace0\rbrace)$.

We are ready to consider a map $\calp$ given by
\begin{equation}\label{calP}
 \calp\left(\begin{array}{l}
                             x_0\\
                             x_1
\end{array}\right)=              \alpha_0x_0 +\int_{-d}^0\bar\mu(d\theta)x_1(\theta):
\end{equation}
this map is well defined as a map $\calp:\cald\rightarrow \R^n$. 
\newline So, given any Borel measurable and bounded function $\bar\phi:\R^n\rightarrow\R$ we can define, $\phi :\cald\rightarrow\R$, by setting
\begin{equation}\label{fiS-gen}
\phi(x)=\bar\phi(\calp(x)) \quad
\forall x=\left( \begin{array}{l}
x_0 \\x_1
\end{array}
 \right)\in \cald,
 \end{equation}
so that $\phi=\bar\phi\circ \calp$.

\noindent The aim now is to prove that the transition semigroup $(R_t)_t$ maps bounded and measurable functions defined by (\ref{fiS-gen}) into differentiable functions: this is the analogous of Theorem \ref{prop-reg} in this more general case of dependence on the past.

\begin{proposition}\label{prop-reg-gen}
Let $(R_t)_{t\geq 0}$ be the transition semigroup  related  to equation (\ref{eq-astr-mild-OU-cald}) defined accordingly to  formula (\ref{ornstein-semS-gen}). 
Let $\phi$ be a bounded and measurable function defined as in (\ref{fiS-gen}) with $\bar\phi \in B_b(\R^n)$.
Then $R_t[\phi]:\cald\rightarrow\R$ is differentiable. Moreover, if \eqref{ass-1} holds, there exists $C>0$ such that for all $h\in \cald_t$ we get the estimate
 \begin{equation}\label{stima_der-gen-Ban}
  \vert\nabla R_t[\phi](x)h\vert \leq C \frac{\vert h\vert}{\sqrt{t}}\Vert \phi\Vert_\infty.
 \end{equation}
\end{proposition}
\dim 
Following the proof of Theorem \ref{prop-reg} we compute, for any $h\in \cald$
\begin{align}\label{Rt-lip-k}
 &R_t[\phi](x+h)- R_t[\phi](x)\\ \nonumber
=&\int_{\R^n}\bar\phi(z+\calp e^{tA}(x+h))\caln(0,\calp Q_t\calp^*)(dz)-\int_{\R^n}\bar\phi(z+\calp e^{tA}x)\caln(0,\calp Q_t\calp^*)(dz)\\ \nonumber
=&\int_{\R^n}\bar\phi(z+\calp e^{tA}x)\caln(\calp e^{tA}h,\calp Q_t\calp ^*)(dz)-
\int_{\R^n}\bar\phi(z+\calp e^{tA}x)\caln(0,\calp Q_t\calp^*)(dz).
\end{align}
 We have to show that the Gaussian measures $\caln(\calp e^{tA}h,\calp Q_t\calp^*)(dz)$
and $\caln(0,\calp Q_t\calp ^*)(dz)$ are equivalent; we will show that under our assumptions the covariance operator $\calp Q_t\calp ^*$ is non degenerated. 
We have to estimate from below the covariance operator: similarly to  \eqref{eq:cov}
we have
\begin{equation*}\<\xi,\calp Q_t \calp^*\xi\>_{\R^n}=\int_0^t \<\xi,\calp e^{sA}GG^*e^{sA^*}\calp^*\xi\>_{\R^n}\,ds=\int_0^t |G^*e^{sA^*}\calp^*\xi|_{\R^n}^2	\,ds=\int_0^t \left[\sup_{ |\eta|=1}\< \calp e^{sA}G\eta,\xi\>_{\R^n}^2\right]ds
\end{equation*}
Fixed $\eta\in \mathbb{R}^n$ with $|\eta|=1$ we have, for $t<d$:
\begin{equation}\label{eq:P-gen} \calp e^{t A}G\eta=\alpha_0 \bar y(t)+\int_{-t}^0 \bar y(t+\theta)\bar\mu(d\theta)
\end{equation}
where $\begin{pmatrix}
\bar y(t)\\
\bar y_t
\end{pmatrix} = e^{tA} \begin{pmatrix}
\sigma \eta\\
0
\end{pmatrix}$ or, by a straight-forward application of the variation of constants formula, again for $t<d$
$$\bar y(t)=e^{t a_0}(\sigma \eta)+\int_0^t e^{(t-s)a_0} \int_{-s}^0 \bar y(s+\theta)a_1(d\theta) \,ds $$
By standard estimates and recalling also that $\bar y(\cdot)$ is continuous we get
\begin{equation}
\label{stima-bar y-gen}
\bar{y}(t)-e^{t a_0}(\sigma \eta)=r_0(t)\; \hbox{ with }\;|r_0(t)|\leq ct a_1([-t,0])
\end{equation}
Notice that by \eqref{ip-eta-stima}, $a_1([-t,0])\rightarrow 0$ as $t\rightarrow 0$, and from \eqref{stima-bar y} we can deduce that 
$\sup_{s\in[0,t]}\vert\bar{y}(s)-e^{s a_0}(\sigma \eta)\vert\leq r_0(t)=o(t)$
Thus
\begin{align*}
\calp e^{t A}G\eta&=\alpha_0 [e^{t a_0}(\sigma \eta)+r_0(t)]+\int_{-t}^0 \left(e^{(t+\theta) a_0}(\sigma \eta)+r_0(t+\theta)\right)\bar \mu(d\theta)\\&=\alpha_0 e^{t a_0}(\sigma \eta)+\int_{-t}^0 e^{(t+\theta) a_0}(\sigma \eta)\bar \mu(d\theta)+r_1(t)
\end{align*}
where $r_1(t)=o(t)$ as $ \rightarrow 0 $.
Thus under \eqref{ass-1}:
$$\calp e^{t A}G\eta=\alpha_0 \sigma \eta+r_2(t) \, \hbox{ where  }\; r_2(t)\sim t \text{ as }t\searrow 0$$
Taking into account \eqref{eq:cov}, \eqref{eq:P-gen} and invertibility of $\alpha_0 \sigma $ we have, for $t\in (0,\epsilon)$ and $\epsilon$ small enough:
$$\<\eta,\calp Q_t\calp^*\eta\>_{\R^n}\geq c |\eta|^2 t$$
\qed

Analogously to Proposition \ref{prop_barQ}, we have to estimate the following operator, that we denote with the same notation used in \eqref{bar-Qst}, for $0\leq s <t\leq T$,
\begin{equation}\label{bar-Qst-ban}
\bar Q^s_t:=\calp e^{sA}Q_{t-s}e^{sA^*}\calp^*=\int_s^t \calp e^{rA}GG^*e^{rA^*}\calp^*\,dr
\end{equation}
\begin{proposition}\label{prop_barQ-ban}
Let $A,\,G$ and $\calp$ be defined respectively in \eqref{A}, \eqref{G-gen} and \eqref{calP} and let \eqref{ass-1} holds true. Let $\bar Q^s_t$ be defined in \eqref{bar-Qst-ban}, $0\leq s<t\leq d$. Then there exists $\bar t>0$ such that $\forall\, 0\,\leq s \leq t\leq \bar t$
\begin{equation}\label{eq:stima_barQ^s_t}
\<\xi,\bar Q^s_t\xi\>_{\R^n}\geq c |\xi|^2 (t-s) \text{ so that }\vert \left(Q^s_t\right)^{-1/2 }\vert \leq c (t-s)^{-1/2}
\end{equation}
\end{proposition}
\dim The proof is completely analogous to he proof of proposition \ref{prop_barQ}, and we omit it. Notice that $\bar t$ depends only on the coefficients of the problem, $ a_0,
\,\alpha_0,\, \mu$ and $d$.
\qed

\section{From the Ornstein-Uhlenbeck to the perturbed Ornstein-Uhlenbeck transition semigroup}
\label{sec:reg_OUpert}

In this Section we consider a perturbed version of the linear delay equation \eqref{eq-ritS_ritlin}, and after the reformulation in the product spaces $H$ and $\cald$ respectively, we notice that under suitable assumptions on the drift we can prove regularizing properties for the perturbed Ornstein-Uhlenbeck transition semigroups, extending the results obtained in Theorem \ref{prop-reg} and in Proposition \ref{prop-reg-gen} for Ornstein-Uhlenbeck transition semigroups.

We start by introducing the semilinear stochastic delay equation. In a complete probability space $(\Omega, \calf,  \P)$ we consider the following stochastic differential equation in $\R^n$ with delay in the state:
\begin{equation}
\left\{
\begin{array}
[c]{l}%
dy(t)  =a_0y(t)dt+\displaystyle \int_{-d}^0y(t+\theta)a_1(d\theta) dt +b(t,y(t), y_t)dt+\sigma dW_t
,\text{ \ \ \ }t\in[0,T] \\
y(0)=x_0\\
y(\theta)=x_1(\theta), \quad \theta \in [-d,0) \, a.e. ,
\end{array}
\right.  \label{eq-ritS}
\end{equation}
where $a_0$ and $a_1$ are as in equation \eqref{eq-ritS_ritlin}, with $a_1$ satisfying \eqref{ip-eta-stima}. Accordingly to Hypothesis \ref{ass}, from now on we consider $\sigma$ invertible.
We consider first the case of initial trajectory $x_1\in L^2([-d,0],\R^n)$ and in equation \eqref{eq-ritS} we consider a drift $b$ with a special form that we describe in the following Hypothesis and that allows to reformulate equation \eqref{eq-ritS} as an abstract evolution equation in the Hilbert space $H$. More precisely, we assume $b$, and consequently $B$  defined in \eqref{B}, is a function with an integral dependence on the past trajectory given by the operator $P$, similarly to (\ref{fiS}) for $\phi$, as it is precisely stated in the following:
\begin{hypothesis}\label{ip-B}
We assume that $b:[0,T]\times\R^n\times L^{2}([-d,0],\mathbb{R}^n)\rightarrow \R^n$ is defined in terms of $\bar b:[0,T]\times R^n\rightarrow \R^n$, setting for all $(t,x_0,x_1)\in[0,T]\times\R^n\times L^{2}([-d,0],\mathbb{R}^n)$
\begin{equation*}
b(t,x_0,x_1)=\bar b\Big( t, P\Big( \begin{array}{l}x_0\\x_1
\end{array}\Big)\Big).
\end{equation*}
We assume that $\bar b$ is continuous and $\forall \,t>0$ $\bar b(t,\cdot):\R^n\rightarrow \R^n$ is Lipschitz continuous and Gateaux differentiable. 
 The map $B:[0,T]\times H \rightarrow \R^n $ is defined  as
 \begin{equation}\label{B-b}
 B(t,x):=\sigma^{-1}b(t,x_0, x_1),
\end{equation}
and so it can be written in terms of a function $\bar B: [0,T]\times \R^n\rightarrow \R^n$ as 
\begin{equation}\label{B}
B(t,x)=\bar B (t,Px), \; \text{where   } \bar B (t,y):=\bar b(t, y),\, \,x=\left( \begin{array}{l}
x_0\\
x_1
\end{array}
\right)\in H,\;y\in \R^n.
\end{equation} 
Due to the assumption on $\bar b$, it turns out that $\bar B$ is bounded and continuous, moreover for every $t\in [0,T]$, $\bar B(t,\cdot)$ is Lipschitz continuous and differentiable. Thus $B$ is bounded, continuous and for every $t\in [0,T]$, $ B(t,\cdot)$ is Lipschitz continuous and G\^ateaux differentiable in $H$.
\end{hypothesis}
So, with $A$ and $G$ defined respectively as in \eqref{A} and \eqref{G}, and $B$ given by \eqref{B-b} and \eqref{B}, equation \eqref{eq-ritS} can be reformulated as a stochastic evolution equation in $H$:
\begin{equation}
\left\{
\begin{array}[c]{@{}l}%
dX^x(t)=AX^x(t)dt+GB(t,X^x(t))+GdW_{t},\text{ \ \ \ }%
t\in\left[ 0,T\right],\\
X^x(0)=x.
\end{array}
\right. \label{eq-astrS}%
\end{equation}
Taking the integral mild form of (\ref{eq-astrS}) we have
\begin{equation}
X^x_t  =e^{tA}x+\int_0^te^{(t-s)A}B(s,X^x_s)ds+\int_0^te^{(t-s)A}GdW_s
,\text{ \ \ \ }t\in[ 0,T] \\
  \label{eq-astr-mild}%
\end{equation}
In Theorem \ref{prop-reg} we have proved regularizing properties for the Ornstein Uhlenbeck transition semigroup 
 $R_t,\,t\geq 0$. 
 \newline Let us now study analogous regularizing properties for the perturbed Ornstein-Uhlenbeck transition semigroup\
 \begin{equation}\label{ornstein-pert-semS}
 P_t[\phi](x)=\E \phi(X^x(t)),\, \phi\in B_b(H).
\end{equation} 
\begin{remark}\label{general-trans-sem-pert}
In equation (\ref{eq-astrS}) we could consider also an initial time $r\neq 0$, namely in mild formulation we get
\begin{equation}
X^{r,x}_t  =e^{(t-r)A}x+\int_r^te^{(t-s)A}B(s,X^x_s)ds+\int_r^te^{(t-s)A}GdW_s
,\text{ \ \ \ }t\in[ 0,T] \\
  \label{eq-astr-mild-r}%
\end{equation}
In this case the transition semigroup related to the process $X^{r,x}$ is given by
\begin{equation}\label{ornstein-pert-semS-r}
 P_{r,t}[\phi](x)=\E \phi(X^{r,x}(t)),\, \phi\in B_b(H),
\end{equation} 
and then it immediately turns out that
$$
 P_{t}[\phi](x) =P_{0,t}[\phi](x)
$$
We underline that if $B=0$ associated to equation \eqref{eq-astr-mild} we have the Ornstein-Uhlenbeck transition semigroup $R_t,\,t\geq 0$, and because  the solution of equation \eqref{eq-astr-mild-r} is omogeneous in time, the Ornstein Uhlenbeck transution semigroup associated to equation \eqref{eq-astr-mild-r} with $B=0$ can be written as $R_{t-r},\,t\geq r$.
\newline In Sections \ref{sec-mildKolmo} and \ref{sec-control}, when dealing with HJB equations and applications to stochastic optimal control, we need to handle transition semigroups associated to processes with initial time not necessarily given by $0$. It is immediate that the reguarizing properties that we prove in this section can be extended to this case.
\end{remark}
We seek for regularizing properties of the transition semigroup $P_t,\, t\geq 0$ when acting on special functions. The strong assumption is that the drift $B$ and the function for which we are able to prove the regularizing properties have the same dependence on the past. Namely, when $B(t,\cdot)$ is Lipschitz continuous and Gateaux differentiable, it turns out that if $\phi$ is a bounded measurable function built as in (\ref{fiS}), then $P_t[\phi]$, for any $t>0$, is a Gateaux differentiable function and for small time $t$ the derivative blows up.
\newline The results are achieved with techniques similar to those of Theorems 4.1 and 4.2 in \cite{Mas1}, adequated to the special context of delay equations and to the case of regularizing properties for transition semigroups when acting on special functions. 
In the proof we will apply the Girsanov Theorem, see e.g. \cite{DP1}, Theorem 10.14.
We set
\[
 V_t^x=\int_0^t\langle B(s,Z^x(s)),dW_s\rangle-\frac{1}{2}\int_0^t \vert B(s,Z^x(s))\vert^2\,ds.
\]
By the Girsanov Theorem we get that $\forall \phi\in B_b(E)$
\begin{equation}
P_{t}\left[  \phi\right]  \left(  x\right)  =\E\left[
\phi\left(  X^{x}(t)\right)\right]=\mathbb{E}\left[
\phi\left(  Z^{x}(t)\right)  \exp V^{x}(t)\right]  ,
\label{transition orn pert}%
\end{equation}
and so $P_{t}$ can be written in terms of the expectation of a function of the
process $Z^{x}$.
\begin{theorem}
\label{teo feller perturbato}
Let us consider the process $X^x$ solution to equation (\ref{eq-astrS}) and let $(P_t)_{t> 0}$ be the related transition semigroup and assume that either \eqref{ass-1} or \eqref{ass-2} holds true.
Let $\phi$ be a bounded and measurable function defined as in (\ref{fiS}), 
by means of $\bar\phi \in B_b(\R^n)$ and assume that \ref{ip-B} holds.
Then $\forall\,t>0,\;P_t[\phi]:H\rightarrow\R$ is differentiable. Moreover, the gradient of $\nabla P_t[\phi]$ can be estimated depending on the non degeneracy assumption we choose. Namely there exists $C>0$ such that for all $h\in H$
\begin{equation}\label{stima_der-pert}
 \hbox{ if \eqref{ass-1} holds then } \vert\nabla P_t[\phi](x)h\vert \leq C \frac{\vert h\vert}{\sqrt{t}}\Vert \phi\Vert_\infty.
\end{equation}
\begin{equation}\label{stimaGder}
 \hbox{ if \eqref{ass-2} holds then } \vert\nabla P_t[\phi](x)h\vert \leq C \frac{\vert h\vert}{t}\Vert \phi\Vert_\infty.
 \end{equation}
\end{theorem}
\dim We prove that for every $\eta\in H$ and for every $\phi$ defined as in (\ref{fiS}) by means of a function 
$\bar\phi\in C^1_b(R^n)$ ( notice that by definition (\ref{fiS}) this implies that $\phi\in C^1_b(H)$ )
\[
\left|  \nabla P_{t}\left[  \phi\right]  \left(  x\right)  \eta\right|
\leq \frac{C}{\sqrt{t}}\left\|  \phi\right\|  _{\infty}\vert \eta\vert ,\ \ \ 0<t\leq T.
\]
We compute $\langle \nabla P_{t}\left[  \phi\right]  \left(  x\right),\eta\rangle$:
\begin{align*}
& \nabla P_{t}\left[  \phi\right]  \left(  x\right)  \eta \\
&  =\lim_{r\rightarrow0}\frac{\mathbb{E}\left[  \bar\phi\left( P Z(t)^{x+r\eta}\right)  \exp \left(V^{x+r\eta}(t)\right)\right]  -\mathbb{E}\left[\bar\phi\left( P Z^{x}(t)\right)  \exp\left( V^{x}(t)\right)\right]  }{r}\\
&  =\lim_{r\rightarrow0}\frac{\mathbb{E}\left[  \bar\phi\left(P  Z^{x+r\eta}(t)\right)  \left(  \exp\left( V^{x+r\eta}(t)\right)-\exp\left( V^{x}(t)\right)\right)\right]  }{r}\\
&\;+\lim_{r\rightarrow0}\frac{\mathbb{E}\left[  \left(\bar\phi\left( P Z^{x+r\eta}(t)\right)  -\bar\phi\left( P Z^{x}(t)\right)\right)  \exp\left( V^{x}(t)\right)\right]  }{r}\\
&  =\mathbb{E}\left[  \bar\phi\left( P Z^{x}(t)\right)  \exp\left( V^{x}(t)\right)\left(\int_{0}^{t}\left\langle \nabla B\left( s, Z^{x}(s)\right)  e^{sA}\eta,dW_{s}\right\rangle \right. 
 \left.  -\int_{0}^{t}\left\langle B\left(  Z^{x}(s)\right) ,\nabla B\left(s,  Z^{x}(s)\right)  e^{sA}\eta\right\rangle ds\right)  \right] \\
 &\; +\mathbb{E}\left[  \left\langle \mathbb{\nabla}%
\bar\phi\left( P Z^{x}(t)\right)  ,Pe^{tA}\eta\right\rangle \exp V^{x}(t)\right] \\
&  =\mathbb{E}\left[  \bar\phi\left( P X^{x}(t)\right)  \int_{0}^{t}\left\langle \nabla B\left(s,  X^{x}(s)\right)  e^{sA}\eta,dW_{s}\right\rangle\right]  +\mathbb{E}\left\langle \mathbb{\nabla}\bar\phi\left( P X^{x}(t)\right)  ,Pe^{tA}\eta\right\rangle 
\end{align*}
The last equality, for what concerns the first term, follows since under the probability measure $\tilde \P$ such that $\dfrac{d\tilde\P}{d\P}=V^x(t)$ the joint law of $\Big(Z^x(\cdot),\, -\dis\int_0^\cdot B(r,Z^x(r))\,dr+W_\cdot\Big)$ coincides with the joint law of $(X^x(\cdot),\, W_\cdot)$ under $\P$, see also \cite{F}, formula (2.12). Moreover it is immediate to see that
\begin{equation}\label{stima1}
\vert \mathbb{E}\left[  \bar\phi\left( P X^{x}(t)\right)  \int_{0}^{t}\left\langle \nabla B\left(s,  X^{x}(s)\right)  e^{sA}\eta,dW_{s}\right\rangle\right]\vert \leq C \Vert \phi\Vert_\infty\vert \eta\vert
\end{equation}
Also we notice that $\forall \,h\in H$, due to the definition of $B$ in (\ref{B}), we can write
\[
\<\nabla B(s,x),h\>_H=\<\nabla \bar B(s,Px),Ph\>_{\R^n}.
\]
This will be useful when evaluating $\mathbb{E}\left\langle \nabla\bar\phi\left(PX^{x}(t)\right)  ,Pe^{tA}\eta\right\rangle $.
Now let $\left(  \xi(t)\right)_t$\ be an $\R^n$-valued bounded predictable process. We define $X^{\varepsilon,x}(t)$ which is the mild solution to the equation
\begin{equation}\label{X-epsilon}
\left\{
\begin{array}
[c]{l}%
dX^{\varepsilon,x}(t)=AX^{\varepsilon,x}(t)dt+GB\left(t,  X^{\varepsilon,x}(t)\right)dt+G\varepsilon\xi(t)dt+GdW_{t},\text{
\ \ \ }t\in\left[  0,T\right] \\
X^{\varepsilon,x}(0)=x.
\end{array}
\right.
\end{equation}
We define the probability measure $Q_{\varepsilon}$ such that%
\[
\dfrac{dQ_{\varepsilon}}{d\mathbb{P}}
=\rho^{\varepsilon}(T),\;
\text{ where }\;
\rho^{\varepsilon}(t)=\exp\left( - \varepsilon\int_{0}^{t}\xi(\sigma)dW_{\sigma}-\dfrac{\varepsilon^{2}}{2}\int_{0}^{t}\left\vert  \xi(\sigma)\right\vert^{2}d\sigma\right)  .
\]
Since $X$ with respect to $\mathbb{P}$ and $X^{\varepsilon}$ with respect to
$Q_{\varepsilon}$ have the same law, it turns out that
\[
\mathbb{E}\left[\bar\phi\left( P X^{x}(t)\right)\right]  =\mathbb{E}\left[  \bar\phi\left(P
X^{\varepsilon,x}(t)\right)  \rho^{\varepsilon}(t)\right].
\]
By differentiating with respect to $\varepsilon$, at $\varepsilon=0$, and applying the dominated convergence theorem, we get%
\begin{align*}
0  &  =\frac{d}{d\varepsilon}_{\mid\varepsilon=0}\mathbb{E}[\bar\phi\left(
P X^{x}(t)\right)]  =\frac{d}{d\varepsilon}_{\mid\varepsilon=0}\mathbb{E}\left[
\bar\phi\left( P X^{\varepsilon,x}(t)\right)  \rho^{\varepsilon}(t)\right] \\
&  =\mathbb{E}\left\langle \mathbb{\nabla}\bar\phi\left( P X^{x}(t)\right)
,P\overset{\cdot}{X}^{\xi}(t)\right\rangle -\mathbb{E}\left[  \bar\phi\left(P
X^{x}(t)\right)
{\displaystyle\int_{0}^{t}}
\left\langle \xi(\sigma),dW_{\sigma}\right\rangle \right]  ,
\end{align*}
where we have denoted 
\begin{equation}\label{X-dot}
\overset{\cdot}{X}^{\xi}(t):=\dfrac{d}{d\varepsilon}_{\mid\varepsilon=0}X^{\varepsilon,x}(t),\;\mathbb{P}-a.s.. 
\end{equation} 
Having defined $\overset{\cdot}{X}$ as the derivative in $0$ and with respect to $\varepsilon$ of the process $X^\varepsilon$ given in \eqref{X-epsilon}, it turns out that $\overset{\cdot}{X}^{\xi}$ is the unique mild solution to the equation
\begin{equation}
\left\{
\begin{array}
[c]{l}%
d\overset{\cdot}{X}^{\xi}(t)=A\overset{\cdot}{X}^{\xi}(t)dt+G\nabla B\left( t, X(t\right)  \overset{\cdot}{X}^{\xi}(t)d+G\xi(t)dt,\text{ \ \ \ }t\in\left[  0,T\right] \\
\overset{\cdot}{X}^{\xi}(0)=0,
\end{array}
\right.  \label{X punto eq}
\end{equation}
that is $\overset{\cdot}{X}^{\xi}$ solves the integral equation%
\[
\overset{\cdot}{X}^{\xi}(t)=\int_{0}^{t}e^{\left( t-s\right)
A}G\nabla B\left(s,  X(s)\right)  \overset{\cdot}{X}^{\xi}(s)ds+\int_{0}%
^{t}e^{\left(t-s\right)  A}G\xi(s)ds.
\]
By the definition of $A$ and $G$ and by Hypothesis \ref{ip-B} on $B$, it can be easily checked that $\overset
{\cdot}{X}^{\xi}$ is well defined as a process with values in $H$. 
From now on it is fundamental to notice that
\[
\nabla B\left(s,  X(s)\right) \overset{\cdot}{X}^{\xi}(s)=\nabla\bar B\left(s, P X(s)\right)P\overset{\cdot}{X}^{\xi}(s)\>,\]
so that the mild form of equation (\ref{X punto eq}) can be written as
\begin{equation}\label{X punto eq mild}
 \overset{\cdot}{X}^{\xi}(t)=\int_{0}^{t}e^{\left( t-s\right)
A}G\nabla \bar B\left(s,  X(s)\right) P \overset{\cdot}{X}^{\xi}(s)ds+\int_{0}%
^{t}e^{\left( t-s\right)  A}G\xi(s)ds.
\end{equation}
Up to now we have proved that for every bounded and predictable process $\left(\xi(t)\right)  _{t}$
\begin{equation}
\mathbb{E}\left\langle \mathbb{\nabla}\bar\phi\left( P X^{x}(t)\right),\,P\overset{\cdot}{X}^{\xi}(t)\right\rangle =\mathbb{E}\left[ \bar\phi\left(PX^{x}(t)\right){\displaystyle\int_{0}^{t}}\left\langle \xi(\sigma),dW_\sigma\right\rangle \right]  .
\label{X punto derivata}%
\end{equation}
This equality, following the proof of Theorem 4.3, pp. 401-402 in \cite{Mas1}, can be extended by density to predictable $\R^n$-valued processes $\left(  \xi(t)\right)_t$ such that $\mathbb{E}\dis\int _{0}^{T}\left\|  \xi(s)\right\|^2ds$ is finite. 
Now we look for a predictable process $\xi\in L^{2}\left(
\Omega\times\left[  0,T\right]  ,\R^n\right)  $ such that%
\[
P \overset{\cdot}{X}^{\xi}(t)=Pe^{tA}\eta.
\]
Let us consider the deterministic controlled system%
\begin{equation}
\left\{
\begin{array}
[c]{l}%
\dfrac{dz(s)}{ds}=Az(s)+Gu(s),\\
z(0)=0,
\end{array}
\right.  \label{sist det}%
\end{equation}
where $u\in L^{2}\left(  \left[  0,T\right]  ,\Xi\right)  $. The solution of
(\ref{sist det}) is given by%
\begin{equation}
z(s)=\int_{0}^{s}e^{\left(  s-r\right)  A}Gu(r)dr.
\label{sist det sol}
\end{equation}
We claim that for all $u\in u\in L^{2}\left(  \left[  0,T\right]  ,\Xi\right)  $ there exists $\xi$ such that $P\overset{\cdot}{X}^{\xi}(s)=Pz(s)$ for every $s\in\left[  0,t\right]  $. Indeed, let us take
\[
\xi(s)=u(s)-\nabla B\left(s,  X(s)\right)  z(s)
=u(s)-\nabla \bar B\left(s, P X(s)\right) P z(s),
\]
$X$ being the solution of \eqref{eq-astr-mild}.
For such a process $\xi$, by considering the mild form of $\overset{\cdot}{X}$ given by (\ref{X punto eq mild}), we get%
\[
P\overset{\cdot}{X}^{\xi}(s)-Pz(s)=\int_{0}^{s}Pe^{\left(  s-r\right)
A}G\nabla \bar B\left(r,X(r)\right)  \left[P  \overset{\cdot}{X}^{\xi}(r)-Pz(r)\right]  dr.
\]
By Gronwall lemma $P\overset{\cdot}{X}^{\xi}(s)-Pz(s)=0$ for every $s\in\left[  0,t\right]  $. Recall that we are looking for $\xi\in\R^n$ such that $P\overset{\cdot}{X}^{\xi}(t)=Pe^{tA}G\eta$, and now we know that this is enough to prove that there exists $u\in L^2([0,T],\Xi)$ such that $Pz(t)=Pe^{tA}G\eta$. Since in Proposition \ref{prop-reg} we have proved that $\bar Q_t$ is non degenerate, we get that $\forall \, t>0$, $\operatorname{Im} Pe^{tA} \subset \operatorname{Im} \bar Q_{t}^{1/2} $: this inclusion implies that  there exists a control $u\in L^{2}\left(  \left[0,T\right]  ,\R^n\right)  $ such that $Pz(t)=Pe^{tA}\eta$. Indeed $\operatorname{Im} \bar Q_{t}^{1/2} =\operatorname{Im}  P Q_{t}^{1/2} =\operatorname{Im} P\call_t$, where $\call_tu=\displaystyle\int_{0}^{t}e^{\left(  t-r\right)  A}Gu(r)dr,$ ( for the last equality see e.g. \cite{DP1} ). So for such a control $u$, by taking $\xi(s)=u(s)-\nabla \bar B\left(s, P X(s)\right) P z(s)$, $0<s<t$, we get that $P\overset{\cdot}{X}^{\xi}(s)=Pe^{tA}\eta$ and that
\[
\mathbb{E}\left\langle \mathbb{\nabla}\bar\phi\left(P  X^{x}(t)\right)
,e^{tA}\eta\right\rangle =\mathbb{E}\left\langle \mathbb{\nabla}%
\bar\phi\left(  X^{x}(t)\right)  ,P\overset{\cdot}{X}^{\xi}(t)\right\rangle
=\mathbb{E}\left[\bar\phi\left(  X^{x}(t)\right)
\int_{0}^{t}
\left\langle \xi(\sigma),dW_{\sigma}\right\rangle\right] .
\]
Moreover
\begin{align*}
\mathbb{E}\left\vert
{\displaystyle\int_{0}^{t}}
\left\langle \xi(\sigma),dW_{\sigma}\right\rangle \right\vert  &
\leq\mathbb{E}\left(
{\displaystyle\int_{0}^{t}}
\left\vert \xi(\sigma)\right\vert^{2}d\sigma\right)  ^{1/2}\\
 &  \leq\mathbb{E}\left(
 {\displaystyle\int_{0}^{t}}
 \vert u(\sigma)\vert^{2}d\sigma\right)  ^{1/2}%
 +\mathbb{E}\left(
 {\displaystyle\int_{0}^{t}}
 \left\vert \nabla \bar B\left( \sigma,P X(\sigma)\right) P z(\sigma)\right\vert^{2}d\sigma\right)^{1/2}\\
 &\leq C\mathbb{E}\left(
 {\displaystyle\int_{0}^{t}}
 \vert u(\sigma)\vert^{2}d\sigma\right)  ^{1/2}
\end{align*}
So for functions $\phi$ defined in (\ref{fiS})
\[
\left\vert \nabla P_{t}\left[  \phi\right]  \left(  x\right)
\eta\right\vert=\left\vert \nabla P_{t}\left[  \bar\phi\right]  \left(  Px\right)
\eta\right\vert  \leq C\Vert\phi\Vert_\infty 
\left(\int_{0}^{t}\vert u(\sigma)\vert^{2}d\sigma\right)  ^{1/2} .
\]
Since the left hand side does not depend on the control $u$, on the right hand side we can take the infimum over all controls $u$ that in the deterministic linear controlled system (\ref{sist det}) steers the initial state $0$ to $Pe^{tA}\eta$ in time $t$. The energy to steer $0$ to $Pe^{tA}\eta$ in time
$t$ is given by%
\[
\mathcal{E}\left(  t,Pe^{tA}\eta\right)  =\min\left\{  \left(  \int_{0}%
^{t}\vert u(s)\vert^{2}ds\right)  ^{1/2}:z\left(  0\right)
=0,\text{ }z\left(  t\right)  =Pe^{tA}\eta\right\}
\]
and  $\mathcal{E}\left(  t,Pe^{tA}\eta\right) =\left\Vert
\bar Q_{t}^{-1/2}e^{tA}\eta\right\Vert $. 
So for functions $\phi$ defined in (\ref{fiS}) by means of $\bar\phi\in C^1_b(\R^n)$ and for 
$\eta\in H$%
\begin{equation}
\left\vert \nabla P_{t}\left[  \bar\phi\right]  \left( P x\right)
\eta\right\vert \leq C\left\Vert \bar Q_{t}^{-1/2}Pe^{tA}\right\Vert\vert \eta\vert \sup_{x\in \R^n}\left\vert
\bar\phi\left(  x\right)  \right\vert .\label{strong feller per perturbati}%
\end{equation}
We claim now that (\ref{strong feller per perturbati}) can be extended to every $\bar\phi\in
C_{b}\left(\R^n\right)  $. 
\newline By convolutions we can uniformly approximate $\bar \phi\in C_b(\R^n)$ with a sequence of functions $\bar \phi_n\in C^\infty(\R^n)$, uniformly bounded by $\Vert \bar\phi\Vert_\infty$. Setting
$$
\phi_n(x):=\bar\phi_n(Px),\,x\in H
$$
we build a sequence of functions $(\phi_n)_n$ infinitely many times differentiable in $H$.
By (\ref{strong feller per perturbati}), we get for every $x,y\in H$
\begin{align*}
\left\vert P_{t}\left[  \phi_{n}\right]  \left(  x\right)  -P_{t}\left[\phi_{n}\right]  \left(  y\right)  \right\vert  &  \leq C\frac{\vert x-y\vert_H}{t^\gamma}\Vert \phi_n\Vert_\infty,
\end{align*}
where $\gamma=\dfrac{1}{2}$ if \eqref{ass-1} holds and $\gamma=\dfrac{3}{2} $ if \eqref{ass-2} holds. 
\newline Letting $n$ tend to $\infty$ in the left hand side and since for every $x\in H,\,P_{t}\left[  \phi_n\right]  \left(  x\right) \rightarrow P_{t}\left[  \phi\right]  \left(  x\right) $,  we get that
\[
\left\vert P_{t}\left[  \phi\right]  \left(  x\right)  -P_{t}\left[
\phi\right]  \left(  y\right)  \right\vert   \leq C\frac{\vert x-y\vert_H}{t^\gamma}
\Vert \phi\Vert_\infty
\]
from which it can be deduced the strong Feller property for the semigroup $P_{t}$. 
\newline We still have to prove that for every $\phi\in C_{b}\left(H\right)  $ defined by means of $\bar\phi \in C_b(\R^n)$, $P_{t}\left[  \phi\right]  $ is a Gateaux differentiable function on $H$. Let us consider again the sequence of approximating functions $(\phi_n)_n$ that we have considered above. By previous calculations we get%
\begin{equation}
\nabla P_{t}\left[  \phi_{n}\right]  \left(  x\right)  \eta
=\mathbb{E}\left[  \bar\phi_{n}\left(P  X^{x}(t)\right)  \int_{0}^{t}\left\langle \nabla \bar B\left(s,  X^{x}(s)\right) P e^{sA}\eta,dW_{s}%
\right\rangle \right]  +\mathbb{E}\left\langle \bar\phi_{n}\left(P  X^{x}(t)\right)  ,\int_{0}^{t}\xi(\sigma)dW_{\sigma}\right\rangle .\label{Hx continuo}
\end{equation}
We get
\begin{align*}
  \nabla &P_{t}\left[  \bar\phi_{n}\right]  \left(  x\right)  \eta
-\nabla P_{t}\left[  \bar\phi_{k}\right]  \left(  x\right)  \eta\\
&  =\mathbb{E}\left[  \left( \bar \phi_{n}\left(  X^{x}(t)\right)-\bar\phi_{k}\left(  X^{x}(t)\right)  \right)  \int_{0}^{t}\left\langle \nabla\bar B\left(s,P  X^{x}(s)\right) P e^{sA}\eta,dW_{s}\right\rangle \right]\\
&\;+\mathbb{E}\left\langle \left(  \bar\phi_{n}\left( P X^{x}(t)\right)-\bar\phi_{k}\left( P X^{x}(t)\right)  \right)  ,\int_{0}^{t}\xi(\sigma)dW_{\sigma}\right\rangle ,
\end{align*}
and the right hand side tends to $0$ in a ball of radius equal to $1$, uniformly with respect to $\eta\in H$.
So there exists $H^{x}\in H$ such that $\nabla P_{t}\left[  \phi_{n}\right]  \left(  x\right)\rightarrow H^{x}$ as $n\rightarrow\infty.$ By (\ref{Hx continuo}), for every $\eta\in H$, the map $x\mapsto H^{x}\eta$ is continuous as a map from $H$ to $\R$. By the estimate (\ref{strong feller per perturbati}) we get that
\[
\left\vert H^{x}\eta\right\vert \leq C\frac{\vert \eta\vert}{t^\gamma}\Vert \phi\Vert_\infty .
\]
It remains to show that $P_{t}\left[  \phi\right]  $ is Gateaux
differentiable and that $\nabla P_{t}\left[  \phi\right]  \left(
x\right)  \eta=H^{x}\eta$. For every $r>0$ and every $\eta\in\R^n$, we can
write
\[
P_{t}\left[  \phi_{n}\right]  \left(  x+r\eta\right)  -P_{t}\left[\phi_{n}\right]  \left(  x\right)  =\int_{0}^{1}\nabla P_{t}\left[\phi_{n}\right]  \left(  x+rs\eta\right)  r\eta ds.
\]
Letting $n\rightarrow\infty$, we get
\[
P_{t}\left[  \phi\right]  \left(  x+r\eta\right)  -P_{t}\left[
\phi\right]  \left(  x\right)  =\int_{0}^{1}H^{x+rs\eta}r\eta ds.
\]
If we divide both sides by $r$ and we let $r$ tend to $0$, by dominated
convergence and by the continuity of $H^{x}\eta$ with respect to $x$, we see
that $P_{t}\left[  \phi\right]  $ is Gateaux differentiable and that $\nabla P_{t}\left[  \phi\right]  \left(  x\right)  \eta=H^{x}\eta$.
Moreover the following estimate holds true: for every $x\in \R^n$ and $\eta\in H$ there exists a constant $C>0$ such that
\[
\left\vert \nabla P_{t}\left[  \phi\right]  \left(  x\right)
\eta\right\vert \leq C\frac{\vert \eta\vert}{t^\gamma}\Vert \phi\Vert_\infty .
\]
\qed

\medskip

We consider also equation (\ref{eq-ritS}) for where the dependence on the past in the drift $b$ is given in terms of the regular measure $\mu$, that satisfies Hypothesis \ref{ip-mu}: in this case equation (\ref{eq-ritS}) will be reformulated in an abstract way in the Banach space $\cald$, when also the initial conditions belong to $\cald$, and we will be able to consider regularizing properties of the perturbed Ornstein Uhlenbeck transition semigroup when acting on special functions defined  as in \eqref{fiS-gen}.
\newline This time on the drift $b$ of equation (\ref{eq-ritS_ritlin}) we assume the following:
\begin{hypothesis}\label{ip-B-gen}
We assume that $b:[0,T]\times\cald\rightarrow \R^n$ is defined in terms of $\bar b:[0,T]\times R^n\rightarrow \R^n$: for all $\Big(t,\Big(\begin{array}{l}x_0\\x_1\end{array}\Big)\Big)\in[0,T]\times\cald)$
\begin{equation*}
b(t,x_0,x_1)=\bar b\Big( t, \alpha_0x_0+\int_{-d}^0x_1(\theta)\bar\mu(d\theta)\Big)\Big), \text{ with }x=\Big(\begin{array}{l}x_0\\x_1\end{array}\Big);
\end{equation*}
and where $\mu$ satisfies Hypothesis \ref{ip-mu}, and $\bar\mu$ is defined by setting $\forall\, A \in \calb([-d,0])$ $\mu(A)=\alpha_{0}\delta_0(A)+\bar\mu(A\setminus\left\lbrace 0\right\rbrace)$.
\newline We assume that $\bar b$ is continuous and $\forall \,t>0$ $\bar b(t,\cdot):\R^n\rightarrow \R^n$ is Lipschitz continuous and differentiable 
 The map $B:[0,T]\times \cald \rightarrow \R^n $ is defined  as in \eqref{B-b}, and so it can be written in terms of a function $\bar B: [0,T]\times \R^n\rightarrow \R^n$ as 
\begin{equation}\label{-genB}
B(t,x)=\bar B \Big(t,\alpha_0x_0+\int_{-d}^0x_1(\theta)\bar\mu(d\theta)\Big), \quad \bar B \Big(t,y\Big):=\sigma^{-1} \bar b\Big(t, y\Big),\; y\in \R^n.
\end{equation} 
Due to the assumption on $\bar b$, it turns out that $\bar B$ is continuous and bounded, and moreover for every $t\in [0,T]$, $\bar B(t,\cdot)$ is Lipschitz continuous and differentiable, and so $B$ is bounded, continuous and for every $t\in [0,T]$, $ B(t,\cdot)$ is Lipschitz continuous and G\^ateaux differentiable in $\cald$.
\end{hypothesis}

Equation  (\ref{eq-ritS}) can be reformulated formally as a perturbed Ornstein-Uhlenbeck process in $\cald$
\begin{equation}
\left\{
\begin{array}
[c]{l}
dX^x(t)  =AX^x(t) dt+GB(t,X^x(t))+GdW_t
,\text{ \ \ \ }t\in[ 0,T] \\
X^x(0)  =x=\left(\begin{array}{l}x_0\\x_1\end{array}\right),
\end{array}
\right.   \label{eq-astrS-gen}%
\end{equation}
and in integral mild form we have
\begin{equation}
X^x(t)  =e^{tA}x+\int_0^te^{(t-s)A}GB(s,X^x(s))ds+\int_0^te^{(t-s)A}GdW_s
,\text{ \ \ \ }t\in[ 0,T] .\\
  \label{eq-astr-mild-gen}%
\end{equation}
For what concerns the drift $GB :[0,T]\times \cald \rightarrow \cald$, it can be proved that, see e.g. \cite{FZ},
$$
\int_{0}^te^{(t-s)A} G B(s,  Y)\, ds \in \cald\,  \quad \text{ if }\;Y\in \cald.
$$
So the perturbed Ornstein-Uhlenbeck process $X^x$ solution to equation \eqref{eq-astrS-gen} is a process in $\cald$ if the initial condition $x\in \cald$.
\newline  Notice also that in terms of $\calp$ the drift $B$ in \eqref{eq-astrS-gen} can be written as
\begin{equation}\label{B-gen}
B(t,y)=\bar B (t,\calp ),\; y\in\cald:
\end{equation} 
this is similar to formula \eqref{-genB}, the difference is only in the use of $\calp$ instead of $P$.

The associated perturbed Ornstein-Uhlenbeck transition semigroup $P_t, \,t\geq0$, is defined by setting, for every bounded measurable function $f:\cald\rightarrow\R$, and for every $x\in \cald$,
\begin{equation}
 \label{ornstein-semS-pert}
P_t[f](x)=\E [f(X^x(t))].
\end{equation}
We now prove that the transition semiproup $(P_t)_{t>0}$ is regularizing also in $\cald$ when acting on special functions, as the Ornstein Uhlenbeck process $(R_t)_{t>0}$ is regularizing in $\cald$, see Proposition \ref{prop-reg-gen}.
Namely, when $B(t,\cdot)$ is Lipschitz continuous and Gateaux differentiable, it turns out that if $\phi$ is a
bounded measurable function 
like the ones defined in (\ref{fiS-gen}), then for any $t>0$ $P_t[\phi]$ is a Gateaux differentiable function and the following estimates hold true:
 \begin{equation}\label{stima_der_pert-gen}
  \vert\nabla P_t[\phi](x)h\vert \leq C \frac{\vert h\vert}{\sqrt{t}}\Vert \phi\Vert_\infty,\;\forall\, h\,\in H.
   \end{equation}
The proof of such regularizing properties for perturbed Ornstein Uhlenebeck transition semigroup and the proof of the related estimates (\ref{stima_der_pert-gen}) are achieved similarly to what we have done in Theorem \ref{teo feller perturbato}.
\begin{theorem}
\label{teo feller perturbato gen}
Let us consider the process $X^x$ solution to equation (\ref{eq-astr-mild-gen}) and let $P_t,\,t\geq 0$ be its transition semigroup related defined in 
(\ref{ornstein-pert-semS}), with $B$ satisfying \ref{ip-B-gen}. 
Let $\phi$ be a bounded and measurable function defined as in (\ref{fiS}) with 
$\bar\phi \in B_b(\R^n)$.
Then $P_t[\phi]:\cald\rightarrow\R$ is differentiable. Moreover, if \eqref{ass-1} holds, there exists $C>0$ such that for all $h\in \cald$ the following estimate holds true:
 \begin{equation}\label{stima_der-gen}
  \vert\nabla P_t[\phi](x)h\vert \leq C \frac{\vert h\vert}{\sqrt{t}}\Vert \phi\Vert_\infty.
 \end{equation}
\end{theorem}
\dim
The proof follows the lines of the proof of Theorem 
\ref{teo feller perturbato} and we omit it.
\qed
%
%

\section{Mild solution of Kolmogorov equations}
\label{sec-mildKolmo}

In this Section we consider Kolmogorov equations: by using the regularizing properties shown in the previous Sections we are able to prove existence of mild solutions, in a sense that we are going to specify. We are able to solve Kolmogorov equations both in $H$ and in $\cald$: we focus on the case of $\cald$, similar arguments apply to $H$.

\subsection{Kolmogorov equations related to the perturbed Ornstein-Uhlenbeck process}
\label{subsec-Kolmo}

In this Section we briefly present (linear) Kolmogorov equations of the form
\begin{equation}\label{Kolmo-formale0-lin}
  \left\{\begin{array}{l}\dis
-\frac{\partial v(t,x)}{\partial t}=\call_t^p [v(t,\cdot)](x),\qquad t\in [0,T],\,
x\in \cald,\\
\\
\dis v(T,x)=\phi(x).
\end{array}\right.
\end{equation}
where by $\call_t^p$ we formally denote the generator of the perturbed Ornstein-Uhlenbeck semigroup (\ref{ornstein-semS-pert}) and by $\<\cdot,\cdot\>$ we mean, again formally, the duality $\<\cdot,\cdot\>_{\cald,\cald^*}$:
\begin{align}\label{eq:ell-pert}
 \call_t^p[f](x)&=\frac{1}{2} Tr \;GG^*\; \nabla^2f(x)
+ \< Ax,\nabla f(x)\>+\< GB(t,x),\nabla f(x)\>\\ \nonumber
&=\frac{1}{2} Tr \;GG^*\; \nabla^2f(x)
+ \< Ax,\nabla f(x)\>+\< B(t,x),G^*\nabla f(x)\>_{\R^n} \nonumber.
\end{align}
We notice that we can also consider equation \eqref{Kolmo-formale0-lin} with the generator $\call$ of the Ornstein Uhlenbeck semigroup \eqref{ornstein-semS-gen}, that is with $B\equiv 0$ in \eqref{eq:ell-pert}, in the place of $\call^p_t$. All the results we present here apply to that case.
\newline By the Feynman-Kac formula the solution to equation \eqref{Kolmo-formale0-lin} is given by
\begin{equation}
v(t,x) =P_{t,T}[\phi](x),\qquad t\in [0,T],\
x\in \cald,\label{solmildKolmo_sem-0}
\end{equation}
Next we define some spaces of directionally differentiable functions, similarly to the spaces defined e.g. in \cite{GoMa1}.
\begin{definition}\label{df:Gspaces}
Let $I$ be an interval in $\R$ and let $K$ be a Banach space.
\begin{itemize}
\item
We call $C^{1}_{b}(K)$  the space of all functions $f:K\to \R$ which admit continuous and bounded G\^ateaux derivative. Moreover we call $C^{0,1}_b(I\times K)$ the space of continuous functions $f:I\times K\to \R$ belonging to $C_b(I\times K)$ and such that, for every $t\in I$, $f(t,\cdot )\in C^{1}_b(K)$. 
\end{itemize}
 From now on $I=[0,T]$.
\begin{itemize}
 \item
For any $\alpha\in(0,1)$ and $T>0$ 
we denote by $C^{0,1}_{\alpha}([0,T]\times K)$ the space of functions
$ f\in C_b([0,T]\times K)\cap C^{0,1}_b((0,T]\times K)$ such that
the map $(t,x)\mapsto t^{\alpha} \nabla f(t,x)$
belongs to $C_b((0,T]\times K,K^*)$.
The space $C^{0,1}_{\alpha}([0,T]\times K)$ is a Banach space when endowed with the norm
\[
 \left\Vert f\right\Vert _{C^{0,1}_{\alpha}([0,T]\times K)  }=\sup_{(t,x)\in[0,T]\times K}
\vert f(t,x)\vert+
\sup_{(t,x)\in (0,T]\times K}  t^{\alpha }\left\Vert \nabla f(t,x)\right\Vert_{K^*}.
\]
We will also write $\left\Vert f\right\Vert _{C^{0,1}_{\alpha}}$ if no confusion is possible.
\end{itemize}
\end{definition}
We notice that $v(T-\cdot, \cdot)\in C^{0,1}_{{1/2}}\left([0,T]\times \cald\right)$. Notice that differently from the classical strong Feller property, we cannot prove any further regularity result, even in the case of the Ornstein Uhlenbeck semigroup. Indeed, in the case of the Ornstein-Uhlenbeck semigroup
\[
v(t,x) =R_{T-t}[\phi](x)=\E[\bar\phi (\calp Z^x(T-t))],
\]
but we cannot guarantee that the derivative of $v$ can be written as a function depending on $x$ only through $\calp x$, so we cannot guarantee any further smoothing.
We will be able to prove existence of a mild solution to the semilinear Kolmogorov equation (\ref{Kolmo-formale0}) when the final datum $\phi:\cald \rightarrow \R$ 
depends on $x$ only through $\calp x$, where $\calp$ has been defined in (\ref{calP}).
\newline Namely, on $\phi$ we make the following assumptions:
\begin{hypothesis}\label{ip-Kolmo-lBfi}
The functions $\phi $ is of the form $\phi=\bar\phi\circ \calp$ as in (\ref{fiS-gen}).
Moreover we assume that $\bar\phi\in C_b(\R^n)$.
 \end{hypothesis}

Next we state an existence result of a mild solution to equation (\ref{Kolmo-formale0-lin}), following \cite{DP1}.
\begin{theorem}\label{teo:esistenzaKolmo_lin}
Let $X^x$ be defined in \eqref{eq-astr-mild-gen}, let $B$ satisfy Hypothesis \ref{ip-B-gen} and let $\phi$ satisfy Hypothesis \ref{ip-Kolmo-lBfi} 
Then if \eqref{ass-1} hold, then the semilinear Kolmogorov equation (\ref{Kolmo-formale0-lin}) admits a unique mild solution $$v(t,x)=\E \phi (X^{t,x}_T),$$ where $X^{t,x}$ is the perturbed Ornstein Uhlenbeck.
Moreover for every $t\in [0,T)$ $v$ is  differentiable with
\begin{equation}\label{reg-sol-kolmo}
\vert\nabla v(t,x)h\vert \leq C \frac{\vert h\vert}{\sqrt{T-t}}\Vert \phi\Vert_\infty, 
\end{equation}  
that is $v(T-\cdot,\cdot)\in C^{0,1}_{1/2}$.
Finally if the initial datum $\phi$ is also continuously Fr\'echet differentiable,
then $v \in C^{0,1}_{b}([0,T]\times \cald)$ and, for suitable $C_T>0$,
\begin{equation}\label{eq:stimavmainteobis-lin}
\Vert v\Vert_{C^{0,1}_{b}}\le C_T\left(\Vert\phi \Vert_\infty
+\Vert\nabla\phi \Vert_\infty\right)
\end{equation}
\end{theorem}
\dim The proof follows directly from the regularizing properties of the transition semigroup studied in Section \ref{sec:reg_OUpert}, Proposition \ref{teo feller perturbato gen}.

\subsection{Semilinear Kolmogorov equations related to the Ornstein-Uhlenbeck process}
\label{subsec-semKolmo}

In this Section we solve semilinear Kolmogorov equations by adequating the fixed point argument to the case when smoothing properties hold true only when acting on special functions.
\newline With these techniques we are able to consider semilinear Kolmogorov equations which formally are given by
\begin{equation}\label{Kolmo-formale0}
  \left\{\begin{array}{l}\dis
-\frac{\partial v(t,x)}{\partial t}=\call[v(t,\cdot)](x) +
\psi (v(t,x), \nabla v(t,x)G),\qquad t\in [0,T],\,
x\in \cald,\\
\\
\dis v(T,x)=\phi(x).
\end{array}\right.
\end{equation}
By $\call$ we denote the generator of the Ornstein-Uhlenbeck semigroup (\ref{ornstein-semS-gen})
\begin{equation}\label{eq:ell}
 \call[f](x)=\frac{1}{2} Tr \;GG^*\; \nabla^2f(x)
+ \< x,A^*\nabla f(x)\>, 
\end{equation}
where as before by $\<\cdot,\cdot\>$ we mean the duality between $\cald$ and $\cald^*$, that is $\<\cdot,\cdot\>_{\cald,\cald^*}$.
\newline Notice that in \eqref{Kolmo-formale0} the non linear term $\psi $ doesn't depend on $x$: with our techniques we aren't able to allow dependence on $x$. Moreover in the proof of Theorem \ref{teo:esistenzaKolmo} for the existence of a mild solution of the semilinear Kolmogorov equation it will be clear that in order to treat the nonlinear term $\psi$ explicit computations are needed, and these techniques cannot be adequated to the perturbed Ornstein Uhlenbeck transition semigroup. 
%
\newline By applying formally the variation of constants formula, see e.g. \cite{DP3}, the semilinear Kolmogorov equation (\ref{Kolmo-formale0}) can then be rewritten in its ``mild formulation''
\begin{equation}
v(t,x) =R_{T-t}[\phi](x)+\int_t^T R_{s-t}\left[
\psi(v(s,\cdot),\nabla v(s,\cdot)G)\right](x)\; ds,\qquad t\in [0,T],\
x\in \cald,\label{solmildKolmo_sem}
\end{equation}
We use this formula to give the notion of mild
solution for the semilinear Kolmogorov equation (\ref{Kolmo-formale0}).

\begin{definition}\label{defsolmildKolmo}
We say that a function $v:[0,T]\times \cald\rightarrow\mathbb{R}$ is a mild solution of the semilinear Kolmogorov equation (\ref{Kolmo-formale0}) if the following are satisfied:
\begin{enumerate}
\item $v(T-\cdot, \cdot)\in C^{0,1}_{{1/2}}\left([0,T]\times \cald\right)$;
\item  equality (\ref{solmildKolmo}) holds on $[0,T]\times \cald$.
\end{enumerate}
\end{definition}
Notice that the right hand side of (\ref{solmildKolmo}) is well defined if $v\in C^{0,1}_{{1/2}}\left([0,T]\times \cald\right)$.

On the final datum $\phi$ we assume that Hypothesis \ref{ip-Kolmo-lBfi} holds true, and on the nonlinear term $\psi$ we make the following assumpions:
 \begin{hypothesis}\label{ip-Kolmo-psi}
The function $\psi:\R\times \R^n\rightarrow \R$ is  Lipschitz continuous: there exists a cnstant $L>0$ such that
  \[
 \vert \psi(\xi, h)-\psi(\eta, k)\vert \leq L\left(\vert \xi-\eta\vert+\vert h-k\vert  \right)
  \]
 \end{hypothesis}
Since the transition semigroup $R_t$ is not even strongly Feller we cannot study the existence and uniqueness of a mild solution of equation (\ref{Kolmo-formale0}) as it is done e.g. in \cite{G1} and in \cite{Mas0}.
We then use the partial smoothing property studied in Sections \ref{sec-delayS} and \ref{sec:reg_OUpert}.
\newline The right space where to seek a mild solution will be the space $\Sigma^{1}_{T,{1/2}}\subset C^{0,1}_{\frac{1}{2}}([0,T]\times \cald)$, that we define below. Indeed our existence and uniqueness result will be proved by a fixed point argument in such a space.
An analogous of this space has been introduced in \cite{GoMa1}, where this space has been introduced related to directional derivatives. 
\begin{definition}\label{df:Sigma}
Let $T>0$ and let $\calp$ be defined as in (\ref{calP}).
A function $g\in C_b([0,T]\times \cald)$belongs to $\Sigma^{1}_{T,\frac{1}{2}}$ if there exists a function $\bar g$ defined in $[0,T]\times \R^n$ such that 
$$g(t,x)=\bar g\left(t,\calp e^{tA}x\right),
\qquad \forall (t,x) \in [0,T]\times \cald,
$$
and if, for any $t\in(0,T]$ the function $g(t,\cdot)$ is Fr\'echet differentiable and if there exists a function $\overline{\nabla g}\in C_b((0,T]\times \R^n)$ such that
$$
t^\frac{1}{2} \nabla g(t,x)=\overline{\nabla g}\left(t,\calp e^{tA}x)\right),
\qquad \forall (t,x) \in [0,T]\times \cald.
$$
\end{definition}
It turns out that $\Sigma^{1}_{T,\frac{1}{2}}$ is a closed subspace of $C^{0,1}_{\frac{1}{2}}([0,T]\times \cald)$ endowed with the norm
\[
 \left\Vert f\right\Vert _{C^{0,1}_{\frac{1}{2}}([0,T]\times \cald)  }=\sup_{(t,x)\in[0,T]\times \cald}
\vert f(t,x)\vert+
\sup_{(t,x)\in (0,T]\times \cald}  t^{\frac{1}{2} }\left\Vert \nabla f(t,x)\right\Vert_{\cald^*},
\]
for more details see \cite{GoMa1}.
\newline From what we have proved in Section \ref{sec:reg_OUpert}, it follows that if $\phi=\bar\phi\circ\calp$ with $\bar\phi\in C_b(\R^n)$, then $R_t[\phi] \in \Sigma^{1}_{T,\frac{1}{2}}.$ We state here an analogous result for convolution type terms that appear below in the fixed point argument.
\begin{lemma}\label{lemma_convoluzione}
Let $(R_t)_{t>0}$ be the Ornstein-Uhlenbeck transition semigroup in $\cald$, defined in (\ref{ornstein-semS-gen}), 
let Hypotheses \ref{ip-Kolmo-psi}, \ref{ip-Kolmo-lBfi} and \eqref{ass-1} hold true and let $\tau\leq \min(d,\bar t, T,1)$, where $\bar t$ is assigned in proposition \ref{prop_barQ-ban}. Then for every $g \in \Sigma^{1,G}_{\tau,\frac{1}{2}}$, if we define the function
$\Gamma g:[0,\tau]\times \cald \rightarrow \R$ by
\begin{equation}\label{iterata-primag}
\Gamma g(t,x) =\int_{0}^{t}
R_{t-s} [\psi(g(s,\cdot),\nabla g(s,\cdot)G)](x)  ds,
\end{equation}
then $\Gamma g$ belongs to $\Sigma^{1}_{\tau,\frac{1}{2}}$.
 Hence, in particular, $\Gamma g(t,\cdot)$ is Fr\'echet differentiable for every $t\in (0,\tau]$ and, for all $x\in \cald$,
\begin{equation}\label{stimaiterata-primag}
\left\vert \nabla (\Gamma g(t,\cdot))(x) \right\vert   \leq
    C\left(t^{\frac{1}{2}}+\Vert g \Vert_{C^{0,1}_{\frac{1}{2}}}\right).
 \end{equation}
\end{lemma}
 \dim
We start by proving that if $g\in \Sigma^{1}_{\tau,\frac{1}{2}}$, then $\Gamma g\in \Sigma^{1}_{\tau,\frac{1}{2}}$.
We have
\begin{align*}
  \int_{0}^{t}
R_{t-s} \left[\psi\left(g(s,\cdot),\nabla g(s,\cdot)G\right)\right](x)  ds
&=\int_{0}^{t}
\int_H \psi\left(g(s,z),\nabla g(s,z+e^{(t-s)A}x)G\right)
\caln(0,Q_{t-s})(dz)
\end{align*}
By the definition of $\Sigma^{1}_{T,\frac{1}{2}}$, there exist $\bar g$ and $\overline{\nabla g}$ such that
$$
g(s,x)=\bar g\left(s, \calp e^{sA}x\right), \quad s^{\frac{1}{2}} \nabla g(s,x) = \overline{\nabla g}\left(s, \calp e^{sA}x\right),
$$
and so
\begin{align}\label{eq:nablaCshiftg}
&g(s,z+e^{(t-s)A}x)=\bar g\left(s, \calp e^{sA}z+\calp e^{tA}x\right), \nonumber\\ 
 &s^{\frac{1}{2}} \nabla g(s,z+e^{(t-s)A}x)
 =\overline{\nabla g}\left(s, \calp e^{sA}z+\calp e^{tA}x\right)
 \qquad \forall t\ge s>0, \; \forall x,z \in H.
\end{align}
Hence the function $\overline{\Gamma g}$ associated to $\Gamma g$ is
\begin{align*}
\overline{\Gamma g} (t,y)&=\int_{0}^{t}
\int_{H} \psi\left( \bar g\left(s, \calp e^{sA}z+y\right),s^{-\frac{1}{2}}
\overline{\nabla g}\left(s, \calp e^{sA}z+y\right)
\right)\caln(0,Q_{t-s})(dz)
\end{align*}
and, by Lipschitz assumptions on $\psi$,
$$
\Vert\overline{\Gamma g}\Vert_\infty \le C\int_{0}^{t}\left( 1 +(1+
s^{-\frac{1}{2}} ) \Vert g \Vert_{C^{0,1}_{\frac{1}{2}}}\right)ds
$$
We compute the derivative:
\begin{align}\label{limite_der}
&\<\nabla \Gamma g, k\>=\< \nabla \int_{0}^{t} R_{t-s}\left[\psi\left( g(s,\cdot),\nabla g(s,\cdot)G)\right)\right](x)\,ds,k\>\\ \nonumber
 &\lim_{\alpha\rightarrow 0}\dfrac{1}{\alpha}
 \left[\int_{0}^{t}
 R_{t-s}\left[\psi\left( g(s,\cdot),\nabla g(s,\cdot)G\right)\right](x+\alpha k)ds -
\int_{0}^{t}
 R_{t-s}\left[\psi\left(g(s,\cdot), \nabla g(s,\cdot)G)\right)\right](x)  ds\right].
\end{align}
From \eqref{eq:nablaCshiftg},
\begin{align*}
&\int_{0}^{t}
 R_{t-s}\left[\psi\left(g(s,\cdot), \nabla g(s,\cdot)G\right)\right](x+\alpha k)ds \\
&=\int_{0}^{t} \int_{\R^n}
\psi\left(\overline g\left(s,z+\calp e^{tA}x\right),s^{-\frac{1}{2}}
\overline{\nabla g}\left(s,z+\calp e^{tA}x\right)G
\right)
\caln\left(
\calp e^{(t-s)A}\alpha k,\calp e^{sA}Q_{t-s}e^{sA*}\calp ^*\right)(dz)ds
\\[2mm]
\end{align*}
Arguing in a similar way on the second term in \eqref{limite_der}, we get that the derivative in \eqref{limite_der} can be rewritten as
\begin{align}\label{limite_der_1}
& \<\nabla  \int_{0}^{t} R_{t-s}\left[\psi\left(\overline{ g}(s,\cdot),\overline{\nabla g}(s,\cdot)G\right)\right](x)\,ds,k\>\\
 &=\lim_{\alpha\rightarrow 0}\dfrac{1}{\alpha}
 \int_{0}^{t}\left[\int_{\R^n}
\psi\left(\overline g\left(s,z+\calp e^{tA}x\right),s^{-\frac{1}{2}}
\overline{\nabla g}\left(s,z+\calp e^{tA}x\right)G\right)
\caln\left(
\calp e^{(t-s)A}\alpha k,\calp e^{sA}Q_{t-s}e^{sA*}\calp ^*\right)(dz)\right.\nonumber \\\nonumber
&\left. -\int_{\R^n}\psi\left(\overline g\left(s,z+\calp e^{tA}x\right),s^{-\frac{1}{2}}
\overline{\nabla g}\left(s,z+\calp e^{tA}x\right)G\right)\caln\left(
\calp e^{(t-s)A}\alpha k,\calp e^{sA}Q_{t-s}e^{sA*}\calp ^*\right)(dz)
\right]\,ds.
\end{align}
As already done in \eqref{bar-Qst-ban}, we set
\[
\bar Q _t^s :=\calp e^{sA}Q_{t-s}e^{sA^*}\calp^*;
\]
with this notation the Gaussian measures in \eqref{limite_der_1} are equivalent if and only if 
\begin{equation}
\operatorname{Im}\calp e^{(t-s)A}  \subset \operatorname{Im}\bar Q _t^s.
\end{equation}
This is true for $t\leq \min (d, \bar t)$, where $\bar t$ is given in Proposition \ref{prop_barQ-ban} and depends only on the coefficients of the problem $\alpha_0,\,a_0,\, \mu$ and $d$. So, setting
\begin{align}
d( t,s,k,z)   &  =\frac{d\caln\left(\calp e^{(t-s) A}k,\bar Q_{t}^s\right)  }{d\mathcal{N}\left(  0,\bar Q_{t}^s\right)  }(z)
\nonumber\\[2mm]
 &  =\exp\left\{  \left\langle (\bar Q_{t}^s)^{-1/2} \calp e^{(t-s)A}k,(\bar Q_{t}^s)^{-{1/2}}z\right\rangle_{\R^n}-\frac{1}{2}\left|  (\bar Q_{t}^s)^{-{1/2}}\calp e^{(t-s)A}k\right|_H  ^{2}\right\}  ,
 \label{eq:density1bisg}
\end{align}
the derivative in \eqref{limite_der_1} can be rewritten as
\begin{align}\label{eq:derCConv}
&\< \nabla  \int_{0}^{t} R_{t-s}\left[\psi\left( g(s,\cdot),\nabla g(s,\cdot)G)\right)\right](x)\,ds,k\>\\[2mm]\nonumber
& =\lim_{\alpha \rightarrow 0}\int_{0}^{t} \int_{\R^n}
\psi\left(\overline{ g}\left(s,z+\calp e^{tA}x\right),s^{-\frac{1}{2}}
\overline{\nabla g}\left(s,z+\calp e^{tA}x\right)G
\right)\frac{d(t,s,\alpha k,z)-1}{\alpha}\\
&\qquad\caln\left(0,\calp e^{sA}Q_{t-s}e^{sA*}\calp ^*\right)(dz)ds
\\[2mm]\nonumber
& =\int_{0}^{t} \int_{\R^n}
\psi\left(\overline{ g}\left(s,z+\calp e^{tA}x\right),s^{-\frac{1}{2}}
\overline{\nabla g}\left(s,z+\calp e^{sA}x\right)G
\right) \left\langle (\bar Q_{t}^s)^{-1/2} \calp e^{(t-s)A}k,(\bar Q_{t}^s)^{-{1/2}}z\right\rangle_{\R^n}\\[2mm]\nonumber
&\qquad\caln\left(0,\calp e^{sA}Q_{t-s}e^{sA*}\calp ^*\right)(dz)ds
\end{align}
Notice that the above formula has the structure required in Definition \ref{df:Sigma}.
\newline In order to prove estimate (\ref{stimaiterata-primag}),
first we set
\begin{align*}
\< \overline{\Gamma g}(t,y),k\>=t^{1/2}
\int_{0}^{t}\int_{\R^n} &\psi \left( \overline g\left(s, z+\calp e^{tA}x\right), s^{-{\frac{1}{2}}} \overline{\nabla g}\left(s, z+\calp e^{tA}x\right) \right)     \left\langle (\bar Q_{t}^s)^{-1/2} \calp e^{(t-s)A}k,(\bar Q_{t}^s)^{-{1/2}}z\right\rangle_{\R^n}\\
&\caln\left(0
\calp e^{sA}Q_{t-s}e^{sA*}\calp ^*\right)(dz)ds;
\end{align*}
 we use the above representation
\begin{align*}
&\left\vert\< \overline{\Gamma g}(t,y),k\>\right\vert
\\[3mm]
 &\leq Ct^{1/2}\int_{0}^{t}
  \int_{\R^n}  \left(1+\left\vert \bar g\left(s, z+y\right)
 \right\vert+\left\vert s^{-{\frac{1}{2}}} \overline{\Gamma g}\left(s, z+y\right) \right\vert\right)
  \left \vert\left\langle (\bar Q_{t}^s)^{-1/2} \calp e^{(t-s)A}k,(\bar Q_{t}^s)^{-{1/2}}z\right\rangle_{\R^n}\right\vert\\[2mm]
  &\quad \caln(0,Q_{t-s})(dz) ds\\
 &\leq Ct^{1/2}\int_{0}^{t}\left(1+s^{-\frac{1}{2}} \left\Vert g \right\Vert_{\Sigma^{1}_{\tau,\frac{1}{2}}}
 \right)
   (t-s)^{-\frac{1}{2}}\vert k\vert_{\cald} \,ds \leq C\left(t^{\frac{3}{2}}+t\left\Vert g \right\Vert_{\Sigma^{1}_{\tau,\frac{1}{2}}} \right)\vert k\vert_{\cald} 
   \leq C \left( 1+\left\Vert g \right\Vert_{\Sigma^{1}_{\tau,\frac{1}{2}}}  \right)\vert k\vert_{\cald} .
\end{align*}
Observe that in the last step we have used estimate
\eqref{eq:stima_barQ^s_t}
\qed

We are now ready to prove our main result of this section on the existence of a mild solution to the semilinear Kolmogorov equation (\ref{Kolmo-formale0}).

\begin{theorem}\label{teo:esistenzaKolmo}
Let $X^x$ be defined in \eqref{eq-astr-mild-gen}, let $B$ satisfy Hypothesis \ref{ip-B-gen} and let $\phi$ satisfy Hypothesis \ref{ip-Kolmo-lBfi} and $\psi$ satisfy Hypothesis \ref{ip-Kolmo-psi}.
Then if \eqref{ass-1} holds, the semilinear Kolmogorov equation (\ref{Kolmo-formale0}) admits a mild solution $v$ according to Definition \ref{defsolmildKolmo}.
Moreover $v$ is unique among the functions $w$ such that $w(T-\cdot,\cdot)\in\Sigma^{1}_{T,1/2}$ and it satisfies, for suitable $C_T>0$, the estimate
\begin{equation}\label{eq:stimavmainteo}
\Vert v(T-\cdot,\cdot)\Vert_{C^{0,1}_{{1/2}}}\le C_T \Vert\bar\phi \Vert_\infty.
\end{equation}
Finally if the initial datum $\phi$ is also continuously 
or Fr\'echet differentiable, then $v \in C^{0,1}_{b}([0,T]\times \cald)$ and, for suitable $C_T>0$,
\begin{equation}\label{eq:stimavmainteobis}
\Vert v\Vert_{C^{0,1}_{b}}\le C_T\left(\Vert\phi \Vert_\infty
+\Vert\nabla\phi \Vert_\infty\right).
\end{equation}
\end{theorem}
\dim First we rewrite (\ref{solmildKolmo_sem}) in a forward way. Namely if $v$ satisfies (\ref{solmildKolmo_sem}) then, setting $w(t,x):=v(T-t,x)$ for any $(t,x)\in[0,T]\times \cald$, we get that $w$ satisfies
\begin{equation}
 w(t,x) =R_{t}[\phi](x)+\int_0^t R_{t-s}[\psi(w(s,\cdot),\nabla w(s,\cdot)G)](x)\; ds,\quad t\in [0,T],\,x\in \cald,
\label{solmildKolmo-forward}
\end{equation}
which is the mild form of the forward HJB equation
\begin{equation}\label{Kolmoformaleforward}
  \left\{\begin{array}{l}\dis
\frac{\partial w(t,x)}{\partial t}=\call [w(t,\cdot)](x) +
\psi (w(t,x),\nabla w(t,x)G),\qquad t\in [0,T],\,
x\in \cald,\\
\\
\dis w(0,x)=\phi(x).
\end{array}\right.
\end{equation}
We use a fixed point argument in $\Sigma^{1}_{T_0,{1/2}}$, with $T_0<\tau$ and $\tau\leq \min(d,\bar t, T,1)$ defined in Lemma \ref{lemma_convoluzione}. $\Sigma^{1}_{T_0,{1/2}}$ is a closed subspace of $C^{0,1}_{{1/2}}([0,T]\times \cald)$, see \cite{GoMa1} where such a space has been studied for directional derivatives.
\newline We define the map $\calc$ on $\Sigma^{1}_{T_0,{1/2}}$ by setting, for $g\in \Sigma^{1}_{T_0,{1/2}}$,
\begin{align}\label{mappaC}
 \calc(g)(t,x):&=R_{t}[\phi](x)+\int_0^t R_{t-s}[
\psi(g(s,\cdot),
\nabla g(s,\cdot)G)\; ds
](x)\\ \nonumber
&=R_{t}[\phi](x)+\Gamma g(t,x),\qquad t\in [0,T_0],\,
x\in \cald,
\end{align}
where in the last equality we have used the same notation of Lemma \ref{lemma_convoluzione}. By \eqref{mappaC}, we will use the smoothing properties of the transition semigroup $(R_t)_t$, also for what concerns the integral term as proved in Lemma \ref{lemma_convoluzione}.
By Proposition \ref{prop-reg-gen} and Lemma \ref{lemma_convoluzione} we deduce that $\calc$ is well defined in $\Sigma^{1}_{T_0,{1/2}}$ and takes its values in $\Sigma^{1}_{T_0,{1/2}}$. Since $\Sigma^{1}_{T_0,{1/2}}$ is a closed subspace  of $C^{0,1}_{{1/2}}([0,T_0]\times\cald)$, once we have proved that $\calc$ is a contraction, by the Contraction Mapping Principle there exists a unique (in $\Sigma^{1}_{T_0,{1/2}}$) fixed point of the map $\calc$, which gives a mild solution  of (\ref{Kolmo-formale0}) in $[0,T_0]$. The same procedure an be iterated in the intervals $[T_0,2T_0],\,[2T_0,3T_0],...$ up to cover the whole interval $[0,T]$.

\noindent Let $g_1,g_2 \in \Sigma^{1}_{T_0,{1/2}}$. We evaluate $\Vert \calc(g_1)-\calc (g_2)\Vert_{\Sigma^{1}_{T_0,{1/2}}}=\Vert \calc(g_1)-\calc (g_2)\Vert_{C^{0,1}_{{1/2}}}$. First of all, arguing as in the proof of Lemma \ref{lemma_convoluzione} we have, for every $(t,x)\in [0,T_0]\times \cald$,
\begin{align*}
  &\vert \calc (g_1)(t,x)- \calc(g_2)(t,x) \vert\\
  & =\left\vert \int_0^t R_{t-s}\left[\psi\left(g_1(s,\cdot),\nabla g_1(s,\cdot)G\right)
 -\psi\left(g_2(s,\cdot),\nabla g_2(s,\cdot)G\right)\right](x)ds\right\vert\\
 &\le \int_0^t Ls^{-{1/2}}\sup_{y \in \cald}\left(\vert s^{{1/2}} (g_1-g_2)(s,y)\vert +s^{{1/2}}\nabla (g_1-g_2)(s,y)G\vert\right) ds\\&
  \leq 2 L t^{{1/2}}\Vert g_1-g_2 \Vert_{C^{0,1}_{{1/2}}}.
\end{align*}
Similarly due to the smoothing properties of $(R_t)_{t>0}$ and arguing similarly to \cite{GoMa1},
\begin{align*}
t^{{1/2}}&\vert \nabla\calc (g_1)(t,x) - \nabla\calc(g_2)(t,x) \vert\\& =
t^{{1/2}}\left\vert \nabla\int_0^t  R_{t-s}\left[\psi\left(g_1(s,\cdot),\nabla g_1(s,\cdot)G\right)-\psi\left(g_2(s,\cdot),\nabla g_2(s,\cdot)G\right)\right](x)ds\right\vert\\
& \leq t^{{1/2}}  L \Vert g_1-g_2 \Vert_{C^{0,1}_{{1/2}}} \int_0^t (t-s)^{-{1/2}}s^{-{1/2}} ds
\le t^{{1/2}}L \beta\left({1/2} , {1/2}\right)\Vert g_1-g_2 \Vert_{C^{0,1}_{{1/2}}}.
\end{align*}
Hence we get
\begin{equation}\label{stima-contr}
 \left\Vert \calc (g_1)-\calc(g_2)\right\Vert _{C^{0,1}_{{1/2}}
  }\leq C\left\Vert g_1-g_2\right\Vert _{C^{0,1}_{{1/2}}
}
\end{equation}
with $C<1$. So the map $\calc$ is a contraction in $\Sigma^{1}_{T_0,{1/2}}$. If we denote by $w$ the unique fixed point of the contraction, then $v:=w(T_0-\cdot,\cdot)$ turns out to be a mild solution of the semilinear Kolmogorov equation (\ref{Kolmo-formale0}) in $[T-T_0
\wedge \delta, T]$, according to Definition \ref{solmildKolmo}.
\newline Since the constant $L$ is independent of $t$, the case of generic $T>0$ follows by dividing the interval $[0,T]$ into a finite number of subintervals of length $T_0\wedge\delta$ sufficiently small.
\newline The estimate (\ref{eq:stimavmainteo}) follows from \cite{GoMa1}. 
\newline Finally the proof of the last statement follows observing that if $\phi$ is continuously Fr\'echet (or Fr\'echet) differentiable, then $R_t[\phi]$ is continuously Fr\'echet differentiable with $\nabla R_t[\phi]$ bounded in $[0,T]\times \cald$.
This allows to perform the fixed point, exactly as done in the first part of the proof,
in $C^{0,1}_b([0,T]\times \cald)$ and to prove estimate (\ref{eq:stimavmainteobis}).
\qed

\section{Application to control}\label{sec-control}

In this Section we apply the previous results to stochastic optimal control problems. We present stochastic optimal control problems that can be refomulated in $\cald$, similarly results apply to problems that can be reformulated in $H$.

Namely we apply the results of Section \ref{sec-mildKolmo} to Hamilton Jacobi Bellman (HJB in the following) equations related to control problems that we are going to present. We can consider  optimal control problems where the state equation is a controlled version of equation \eqref{eq-ritS} and the cost functional reduces to final cost functional: in this case the HJB equation related is a linear Kolmogorov equation with the structure of equation \eqref{Kolmo-formale0-lin}. Moreover we can consider  optimal control problems where the state equation is a controlled version of equation \eqref{eq-ritS_ritlin} and where we can deal with a cost functional with both running and final cost, and where the running cost depends only on the control $u$: in this case the HJB equation related is a semilinear Kolmogorov equation with the structure of equation \eqref{Kolmo-formale0} where in the semilinear term $\psi$ there is no dependence on the state variable $x$.

Let us consider first the following stochastic controlled state equation in $\R^n$ with delay in the state, and with the general dependence on the past that we can treat:
\begin{equation}
\label{eq:SDEx:control-pert}
\left\{
\begin{array}
[c]{l}%
dy(t)  =a_0y(t)dt+\displaystyle \int_{-d}^0y(t+\theta)\eta(d\theta) dt +b(t,y(t), y_t)dt+\sigma u(t)dt+\sigma dW_t
,\text{ \ \ \ }t\in[0,T] \\
y(t_0)=x_0\\
y(\theta)=x_1(\theta), \quad \theta \in [-d,0),
\end{array}
\right. 
\end{equation}
where $b$ satisfies \ref{ip-B-gen}. The control process, $u$ is an $\left(  \mathcal{F}_t\right)  _t$-predictable process with values in a closed and bounded set $U\subset \R^n$. Equation \eqref{eq:SDEx:control-pert} 
can be reformulated in the space $\cald$ as 
\begin{equation}
\label{eq:SDEX:control-pert}\left\lbrace\begin{array}{l}
dX^u(t)=AX^u(t)\,dt +GB(t,X^u(t))+G u(t) dt+G d W(t)\ ,\quad t\in(t_0,T]\\
  X^u(t_0)=\left( \begin{array}{l}x_0\\x_1
  \end{array}\right):=x\ \text{,}\end{array}\right.
\end{equation}
where $B$ and $G$ are defined respectively in \eqref{B-gen} and (\ref{G-gen}).
The solution of this equation will be denoted by $X^{u,t_0,x}$ or simply by $X^{u}$. $X^u$ is also called the state, $T>0,$ $t\in\left[  0,T\right]$ are fixed.
\newline Beside equation (\ref{eq:SDEX:control-pert}), we define the cost
\begin{equation}
J\left(  t_0,x,u\right) =\mathbb{E}\bar\phi\left(\calp X^{u}_T\right). 
\label{cost-semplice}%
\end{equation}
where $\bar\phi:\R\rightarrow\R^n$ is a  bounded and continuous function. We underline the fact that since in equation \eqref{eq:SDEx:control-pert} the drift $b\neq 0$, we can consider only cost functionals where only the final cost is non null, because the HJB equation related is a Kolmogorov equation as \eqref{Kolmo-formale0-lin}: since  the transition semigroup related to equation \eqref{eq:SDEX:control-pert} with $u\equiv 0$ turns out to be a perturbed Ornstein Uhlenbeck treansition semigroup, we are able to consider only linear Kolmogorov equations. We notice that the final cost depends on the trajectory of the state through $\calp$, and this allows more generality than standard final costs, since in some sense the final cost includes a current cost depending only on the state. 
\newline The cost $J$ can be reformulated in an abstract way as
\begin{equation}
J\left(  t,x,u\right)  =\mathbb{E}\phi\left(X^{u}(T)\right),
\label{cost:abstract-pert}%
\end{equation}
where $\phi$ is defined by means of $\bar\phi$ as done in \eqref{fiS-gen}.

Keeping in mind that the HJB equation related to this control problem is the semilinear Kolmogorov equation \eqref{Kolmo-formale0-lin}, the treatment of the control problem goes on in a similar, simpler way to the next case we are going to present.


\medskip

Now we consider the case where in the controlled state equation  \eqref{eq:SDEx:control-pert} the drift $b\equiv 0$, namely we consider
\begin{equation}
\label{eq:SDEx:control}
\left\{
\begin{array}
[c]{l}%
dy(t)  =a_0y(t)dt+\displaystyle \int_{-d}^0y(t+\theta)\eta(d\theta) dt +\sigma u(t)dt+\sigma dW_t
,\text{ \ \ \ }t\in[0,T] \\
y(t_0)=x_0\\
y(\theta)=x_1(\theta), \quad \theta \in [-d,0),
\end{array}
\right. 
\end{equation}
Equation \eqref{eq:SDEx:control} can be reformulated in the space $\cald$ as 
\begin{equation}
\label{eq:SDEX:control}\left\lbrace\begin{array}{l}
dX^u(t)=AX^u(t)\,dt +G u(t) dt+G d W(t)\ ,\quad t\in(t_0,T]\\
  X^u(t_0)=\left( \begin{array}{l}x_0\\x_1
  \end{array}\right):=x\ \text{,}\end{array}\right.
\end{equation}
where $B$ and $G$ are defined respectively in \eqref{B-gen} and (\ref{G-gen}).
\newline In mild formulation the solution of equation (\ref{eq:SDEX:control}) is given by $X^u$ satisfying, for every $t \in [0,T]$,
\begin{equation}
  \label{eq:SDEXmild:control}
  X^u(t)=e^{(t-t_0)A}x+\int_{t_0}^t e^{(t-s)A}G u(s)\,d s+\int_{t_0}^te^{(t-s)A}G \,d W(s).
\end{equation}
Beside equation \eqref{eq:SDEX:control} we consider the cost functional\begin{equation}
J\left(  t,x,u\right)  =\mathbb{E}\int_{t}^{T} g\left(u(s)\right)ds+\mathbb{E}\bar\phi\left(\calp x^{u}_T\right). 
\label{cost}%
\end{equation}
for real functions  $ g$ on $U$ and $\bar\phi$ on $\R^n$.

 Notice that $J$ depends on the state only through the final cost, and the final cost depends on the trajectory of the state, so in some sense it is somehow equivalent to aclassical current cost depending on the state.  

We make the following assumptions on the cost $J$.
\begin{hypothesis}\label{ip costo}
$g:U\rightarrow\R$ is measurable and bounded and $\bar\phi $ is continuous and bounded.
\end{hypothesis}
The cost $J$ can be reformulated in an abstract way as
\begin{equation}
J\left(  t,x,u\right)  =\mathbb{E}\int_{t}^{T}g\left(u(s)\right)ds+\mathbb{E}\phi\left(X^{u}(T)\right),
\label{cost:abstract}%
\end{equation}
where $\phi$ is defined by means of $\bar\phi$ as done in \eqref{fiS-gen}.
\newline The special structure of equation (\ref{eq:SDEX:control})  and of the cost $J$ in \eqref{cost:abstract} leads to a semilinear HJB equation with the structure of the Kolmogorov equation (\ref{Kolmo-formale0}) studied in the previous sections. 
\newline The control problem in strong formulation is to minimize this functional $J$ over all admissible controls $u$.
In the following we denote by $\mathcal{A}_{d}$ the set of admissible controls, that is the $U$-valued predictable processes taking values in $U$.
\noindent We denote by $J^{\ast}\left(  t,x\right)  =\inf_{u\in\mathcal{A}_{d}}J\left(  t,x,u\right)  $ the value function of the problem and, if it exists, we denote by $u^{\ast}$ the optimal control, that is the one realizing the infimum.

We define in a classical way the Hamiltonian function relative to the above problem:%
\begin{equation}\label{def-di-psi}
\psi\left(z\right)  =\inf_{u\in U}\left\{  g\left(u\right)
+zu\right\}\quad \forall z\in \R^d .
\end{equation}
We notice that the Hamiltonian function $\psi$ satisfies Hypothesis \ref{ip-Kolmo-psi}: it is Lippschitz continuous wth respect to its argument.
\newline Moreover we set
\begin{equation}\label{Gamma}
\Upsilon(z)=\{ u:g(u)+zu=\psi(z) \}.
\end{equation}
Notice that $\forall z\, \in\R^d$, $\Gamma(z)$ is not empty since $U$ is compact and $g$ is continuous.
 
We consider the HJB equation related, which has the structure of the semilinear Kolmogorov equation (\ref{Kolmo-formale0}) studied in Section \ref{sec-mildKolmo}, and which is given by 
\begin{equation}\label{HJBformale}
  \left\{\begin{array}{l}\dis
-\frac{\partial v(t,x)}{\partial t}=\call[v(t,\cdot)](x)+\psi (\nabla v(t,x)G),\qquad t\in [0,T],\,
x\in \cald,\\
\\
\dis v(T,x)=\phi(x).
\end{array}\right.
\end{equation}
By applying the variation of constants formula, the HJB equation (\ref{HJBformale}) can be rewritten in its ``mild formulation''
\begin{equation}
v(t,x) =R_{T-t}[\phi](x)+\int_t^T R_{s-t}\left[
\psi(\nabla v(s,\cdot)) G\right](x)\; ds,\qquad t\in [0,T],\
x\in \cald,\label{solmildKolmo}
\end{equation}
The solution $v$ exists and it is unique in $\Sigma^1_{T_0,1/2}$, and we can prove that it admits a semimartingale representation, as stated in the following lemma.
\begin{lemma}\label{lemma:semimart}
Let us consider equation \eqref{eq:SDEx:control} with $a_0\in \operatorname{Mat}(n\times n)$ and $a_1$ an $n\times n$ matrix valued finite regular measure, and let $X^u$ be the solution of its abstract reformulation \eqref{eq:SDEX:control}; let Hypotheses \ref{ip-Kolmo-lBfi}  and \ref{ip costo} hold true. 
If $v$ is the solution to (\ref{solmildKolmo}) and $X^{u} $ is the solution to the controlled stochastic differential equation (\ref{eq:SDEX:control}) with initial time $t_0=t$, there exists a square integrable, $n$-dimensional, adapted process $Z$  such that the process $v\left( s ,X^{u}\left( s\right)  \right) ,\,0\leq t\leq s\leq T $ can be represented as%
\begin{align}
v\left( s,X^{u}\left( s\right)  \right)   &  =v(t,x)+\int_{t}^{s}Z\left(  r\right)  dW\left(  r\right) \label{vsemimartingala}\\
&  +\int_{t}^s\left[  \psi\left(\nabla v\left( X^{u}\left(  r\right)\right) G \right)  +Z\left(r\right)   u\left(  r\right)\right]  dr.\nonumber
\end{align}
Consequently, the process $v\left( s ,X^{u}\left( s\right)  \right), \,0\leq t\leq s\leq T$ is a semimartingale.
\end{lemma}
\dim In the filtered probability space $\left(\Omega,\mathcal{F},\left(\mathcal{F}_s\right)  _{s\geq0}\mathbb{P}\right)  $ we consider equation \eqref{eq:SDEX:control}, where $W$ is an $n$-dimensional standard Brownian motion. The Girsanov theorem ensures that there exists a probability measure $\widetilde{\mathbb{P}}$ such that the process
\[
\widetilde{W}\left(  s\right)  :=W\left( s\right)  +\int_{t}^{s} u\left( r\right)  dr
\]
is an $n$-dimensional standard Brownian motion with respect to $\widetilde{\mathbb{P}}$; we set $\left(\widetilde{\mathcal{F}_{s}}\right) _{s\geq 0}$ the filtration generated by $\widetilde W$ and augmented in the usual way.
\newline Equation \eqref{eq:SDEX:control} with initial time $t_0=t$ can be rewritten as 
\begin{equation}
\left\{
\begin{array}
[c]{l}X^u(s)=AX^u(s)\,ds +Gd \widetilde W (s)\ ,\quad 0\leq t\leq s\leq T\\
  X^u(t)=x,
\end{array}
\right.  \label{eq:SDEX:control1}%
\end{equation}
we also denote its solution by $X^{u,t,x}$.
In the following we write $\widetilde{\mathbb{E}}^{\widetilde\calf_t} [\cdot]=\widetilde{\mathbb{E}} [\cdot|{\widetilde\calf_t}]$. We notice that by \eqref{solmildKolmo}, $\forall \, 0\leq t\leq \tau\leq T$ we can write
\begin{align*}
v(\tau,X^{u,t,x}(\tau))& =\widetilde{\E}\left[\mathbb{\phi}\left(  X^{u,\tau,X^{u,t,x}(\tau)}\left(  T\right)  \right) \right]  \ +\int_\tau^{T}\widetilde{\E}\left[  \psi\left(\nabla v\left(s,X^{u,\tau,X^{u,t,x}(\tau)}\left(  s\right) \right)G\right)  \right]  \, ds\\
&  =\widetilde{\E}^{\widetilde\calf_\tau}\left[\mathbb{\phi}\left(  X^{u,\tau,X^{u,t,x}(\tau)}\left(  T\right)  \right) \right]  \ +\int_\tau^{T}\widetilde{\E}^{\widetilde\calf_\tau}\left[  \psi\left(\nabla v\left(s,X^{u,\tau,X^{u,t,x}(\tau)}\left(  s\right) \right)G\right)   \right]  \, ds\\
&  =\widetilde{\E}^{\widetilde\calf_\tau}\left[\mathbb{\phi}\left(  X^{u,\tau,X^{u,t,x}(\tau)}\left(  T\right)  \right) \right]  \ +\int_t^{T}\widetilde{\E}^{\widetilde\calf_\tau}\left[  \psi\left(\nabla v\left(s,X^{u,\tau,X^{u,t,x}(\tau)}\left(  s\right) \right)G\right)   \right]  \, ds\\
&  -\int_{t}^{\tau}\psi\left(\nabla v\left(s,X^{u,\tau,X^{u,t,x}(\tau)}\left(  s\right) \right)G\right)  ds.
\end{align*}
By the martingale representation theorem
there exists a square integrable process $Z\left(  s\right) ,\,0\leq t\leq s\leq T $ such that the square integrable martingale $$\widetilde{\E}^{\widetilde\calf_\tau}\left[\mathbb{\phi}\left(  X^{u,t,x}\left(  T\right)  \right)+\displaystyle \int_t^{T}\psi\left(\nabla v\left(s,X^{u,\tau,X^{u,t,x}(\tau)}\left(  s\right) \right)G\right)    \, ds\right],\, \tau \in[t,T] $$ can be represented as $v(t,x)+\displaystyle\int_{t}^{\tau}Z\left( r\right)  d\widetilde W\left(r\right),\, \tau \in[t,T]$, and consequently
the process $v\left(  \tau, X^{u,t,x}(\tau) \right)  $ can be written as
\begin{align*}
v\left(  \tau, X^{u,t,x}(\tau)  \right)   &  =v(t,x)+\int_{t}^{\tau}Z\left( r\right)  d\widetilde W\left(r\right)   -\int_{t}^{\tau}\psi\left(\nabla v\left(s,X^{u,\tau,X^{u,t,x}(\tau)}\left(  s\right) \right)G\right)  ds\\
 &  =v(t,x)+\int_{t}^{\tau}Z\left(  r\right)  dW\left(  r\right)  +\int_{t}^{\tau}\left[  Z\left(r\right)   u\left(  r\right)-\psi\left(\nabla v\left(s,X^{u,\tau,X^{u,t,x}(\tau)}\left(  s\right) \right)G\right) \right]  dr
\end{align*}
and so it is a semimartingale in $\left(  \Omega,\mathcal{F},\mathbb{P}%
\right)  $ with respect to the filtration $\left(  \mathcal{F}_{s}\right)
_{s\geq0}$.
\qed 

We have to prove that $v$, solution to \eqref{solmildKolmo}, is the value function of the optimal control problem, and that the optimal control in feedback law is given by means of $\nabla v G$, the derivative in the directions in the image of the operator $G$. These results are collected in the following Theorem: in the proof we take advantage of the fact that equation \eqref{eq:SDEX:control} satisfies the so called structure condition according to which the control affects the system only through the noise. This allows to prove the fundamental relation in terms of backward stochastic differential equations (BSDEs) when the final cost and the Hamiltonian function are differentiable; by an approximation procedure and thanks to the regularizing properties we have proved, we are able to handle a continuous and bounded final cost and Hamiltonian function not necessarily differentiable.
\begin{theorem}
\label{teo-controllo} Let $A$ be defined in \eqref{A}, let $B$ satisfy Hypothesis \ref{ip-B-gen}, let \eqref{ass-1} and \ref{ip costo} hold true. Let $v$ be the mild solution of the HJB equation (\ref{HJBformale}). Then for every $t\in\left[  0,T\right]  $ and $x\in \cald$ and for every admissible control we have 
\begin{equation}\label{fund_rel}
J(  t,x,u )  \geq v\left(  t,x\right)  ,
\end{equation}
where the equality holds if and only if $u$ is optimal. 
Let $\Upsilon_{0} :\cald\times\R^{n}\longrightarrow U$ a Borel measurable selection of $\Upsilon$ defined in \eqref{Gamma} map
If $u$ is an admissible control satisfying%
\begin{equation}
u(s)=\Upsilon_{0}\left(  \nabla v(s X^{u}(s)G\right)  \text{
\ \ \ \ }\mathbb{P}\text{-a.s. for almost every }0\leq t\leq s\leq T ,
\label{u in gammazero}%
\end{equation}
then $J\left(  t,x.u  \right)  =v\left(  t,x\right)  .$
\newline Moreover  the closed loop equation
\begin{equation}
\left\{
\begin{array}
[c]{l}%
dX^{u}(s)=AX^{u}(s)ds+G\Upsilon_{0}\left( \nabla v\left(  s,X^{u}(s)\right)  G \right)  ds+GdW_s,\text{\ \ }0\leq t\leq s\leq T ,\\
X(t)^{u}=x.
\end{array}
\right.  \label{closed loop eq}%
\end{equation}
admits a weak solution in mild sense for a.a. $s \in [t, T)$.
\end{theorem}
\dim We start by smoothing the coefficients of our problem.
 We denote by $(\bar \phi_n)_n$ a sequence of bounded and smooth ( i.e. infinitely many times differentiable) functions $\R^n\rightarrow \R$ such that$$\Vert \bar \phi_n -\bar\phi\Vert_\infty\rightarrow 0 \text{ as }n\rightarrow \infty,\quad\Vert \bar\phi _n\Vert_\infty \leq \Vert \bar\phi \Vert_\infty\, \forall n\geq1$$
Such functions can be built in a standard way by the convolution of $\bar\phi$ with a sequence of smooth kernels. 
Setting $\phi_n:=\bar\phi_n\circ\calp$, we obtain a sequence of functions $\phi_n:\cald\rightarrow\R$ such that
\[
\Vert \phi _n-\phi \Vert_\infty\rightarrow 0,\;\Vert\phi_n\Vert_\infty \leq \Vert \phi \Vert_\infty.
\]
Moreover the functions $\phi_n$ are (infinitely many times) Frechet differentiable. 
\newline In the same way we have approximated $\bar \phi $, we can also approximate $\psi:\R^n\rightarrow \R$ by  making the convolution with smooth kernels and building a sequence of functions $(\psi_n)_n$, infinitely many times differentiable and such that the first order derivative is bounded by the Lipschitz constant $ L$ of $\psi$ uniformly with respect to $n$.
\newline Let us consider, for all $n\geq 1$, the sequence of approximating Hamilton Jacobi Bellman equations given by 
\begin{equation}\label{HJBformale-n}
  \left\{\begin{array}{l}\dis
-\frac{\partial v_n(t,x)}{\partial t}=\call[v_n(t,\cdot)](x)+\psi_n (\nabla v_n(t,x)G),\qquad t\in [0,T],\,
x\in \cald,\\
\\
\dis v_n(T,x)=\phi_n(x),
\end{array}\right.
\end{equation}
which in mild formulation can be rewritten as
\begin{equation*}
v_n(t,x) =R_{T-t}[\phi_n](x)+\int_t^T R_{s-t}\left[
\psi_n(\nabla v_n(s,\cdot)G)\right](x)\; ds,\qquad t\in [0,T],\
x\in \cald.\
\end{equation*}
We consider also the mild form of the equation satisfied by $v_n(T-\cdot, \cdot)$:
\begin{equation}
v_n(T-t,x) =R_{t}[\phi_n](x)+\int_0^t R_{t-s}\left[
\psi_n(\nabla v_n(T-s,\cdot)G)\right](x)\; ds,\qquad t\in [0,T],\
x\in \cald.\label{solmildKolmo-n}
\end{equation}
We argue in a similar way to what done in in \cite{Mas}, but with suitable arrangements due to the fact that we only have smoothing on special functions, to prove that 
\begin{equation}\label{approx-v}
\Vert v_n(T-\cdot, \cdot)-v(T-\cdot, \cdot)\Vert_{C^{0,1}_{1/2}}\rightarrow 0 \text{ as }n\rightarrow \infty.
\end{equation}
Coming into the details of the proof of the convergence in \eqref{approx-v}, by \cite{FT2}, theorem 6.2, equation \eqref{solmildKolmo-n} admits a mild solution $v_n(T-\cdot, \cdot)$ which is jointly continuous and for every
$t>0$, $v_{n}(T-t,\cdot)\in C^{1}\left(  \cald\right)  $. We claim
that $v_n(T-\cdot, \cdot)$ converges uniformly to $v $ in $[0, T_0]$, $0< T_0\leq T$. In view of the application of the theorem of contractions depending on a parameter, we define the maps $\calc_n, \,n\geq 1$ on $\Sigma^{1}_{T_0,{1/2}}$  by setting, for $g\in \Sigma^{1}_{T_0,{1/2}}$,
\begin{equation}\label{mappaCn}
 \calc_n(g)(t,x):=R_{t}[\phi_n](x)+\int_0^t R_{t-s}[
\psi_n(\nabla g(s,\cdot)G)
](x)\; ds,\qquad t\in [0,T_0],\,
x\in \cald.
\end{equation}
For every $n\in\N$, $\calc_n$ turns out to be a contraction in the space $\Sigma^{1}_{T_0,{1/2}}$, indeed let $g_1,g_2 \in \Sigma^{1}_{T_0,{1/2}}$. We evaluate $\Vert \calc_n(g_1)-\calc_n (g_2)\Vert_{\Sigma^{1}_{T_0,{1/2}}}=\Vert \calc_n(g_1)-\calc_n (g_2)\Vert_{C^{0,1}_{{1/2}}}$. Arguing as in the proof of Theorem \ref{teo:esistenzaKolmo}, for every $(t,x)\in [0,T_0]\times \cald$,
\begin{align*}
  &\vert \calc_n (g_1)(t,x)- \calc_n(g_2)(t,x) \vert\\
  & =\left\vert \int_0^t R_{t-s}\left[\psi_n\left(\nabla g_1(s,\cdot)G\right)
 -\psi_n\left(\nabla g_2(s,\cdot)G\right)\right](x)ds\right\vert\\
 &\le \int_0^t Ls^{-{1/2}}\sup_{y \in \cald}s^{{1/2}}\nabla (g_1-g_2)(s,y)G\vert ds\\&\leq 2 L t^{{1/2}}\Vert g_1-g_2 \Vert_{C^{0,1}_{{1/2}}},
\end{align*}
and 
\begin{align*}
t^{1/2}&\vert \nabla\calc_n (g_1)(t,x) - \nabla\calc_n(g_2)(t,x) \vert\\& =t^{1/2}\left\vert \nabla\int_0^t  R_{t-s}\left[\psi_n\left(\nabla g_1(s,\cdot)G\right)-\psi_n\left(\nabla g_2(s,\cdot)G\right)\right](x)ds\right\vert\\
& \leq t^{{1/2}}  L \Vert g_1-g_2 \Vert_{C^{0,1}_{{1/2}}} \int_0^t (t-s)^{-{1/2}}s^{-{1/2}} ds\le t^{{1/2}}L \beta\left({1/2} , {1/2}\right)\Vert g_1-g_2 \Vert_{C^{0,1}_{{1/2}}}:
\end{align*}
what is crucial is that $\psi_n,\, n\geq 1$ are Lipschitz continuous, with the same constant $L$, independent on $n$. Hence we get
\begin{equation}\label{stima-contr-n}
 \left\Vert \calc_n (g_1)-\calc_n(g_2)\right\Vert _{C^{0,1}_{{1/2}}
  }\leq C\left\Vert g_1-g_2\right\Vert _{C^{0,1}_{{1/2}}
}
\end{equation}
with $C<1$. So the maps $(\calc_n)_n$ are contractions in $\Sigma^{1}_{T_0,{1/2}}$. 
 Denote by $ \overline{v} _n$ and $ \overline{v} $ the unique fixed point respectively of $\calc_n$ and of $\calc$. By the theorem of contractions depending on a parameter, we conclude that, as $n\rightarrow\infty,$ $  \overline{v}_n\rightarrow\overline{v}  $ in $C^{0,1}_{{1/2}}([0,T_0]\times \cald)$, see \cite{Mas}. This procedure can be repeated, in the interval $[T_0,2T_0\wedge T]$ and then if $2T_0<T$ it can be repeated in $[2T_0, 3T_0\wedge T]$ and so on: in this way we can conclude that $  \overline{v}_n\rightarrow\overline{v}  $ in $C^{0,1}_{{1/2}}([0,T]\times \cald)$.
\newline Thanks to the approximation of $v$ given in (\ref{approx-v}) and to the fact that the functions $v_n,\,n\geq 1,$ are G\^ateaux differentiable functions with derivative bounded by a constant depending on $n$, see \cite{FT2} and also \cite{FMT} for the path dependent case, it is possible to prove that
\[
Z_n(s)=\nabla v_n(s, X^u(s))G. 
\] 
We have to prove that $Z(s)=\nabla v(s, X^u(s))G$.
We consider the local martingales
\begin{align*}
\zeta_{n}\left(  s,X^{u}\left( s\right)  \right)   &  :=v_{n}\left(s,X^{u}\left( s \right)  \right)  -v_{n}(t,x)\\
&  +\int_{t}^{s}\left[  \psi_{n}\left( \nabla v_{n}\left(r,X^{u}\left(  r\right)\right)  G \right)  -\nabla v_{n}\left(  r,X^{u}\left( r \right)  \right)Gu\left(  r\right)  \right]  dr.
\end{align*}
Since since $v_n\to v$ in $\Sigma^1_{T_0,1/2}$, we can deduce that $v_{n}\left(  s,X^{u}\left(  s\right)  \right)  $  converges a.s to $v\left( s,X^{u}\left(  s\right)  \right)  $, uniformly with respect to the time $s$. Moreover, as $n\to\infty$,
\begin{align*}
\int_{t}^{s}&\left[  \psi_{n}\left( \nabla v_{n}\left(r, X^{u}\left(  r\right) \right)  G \right)  -\nabla v_{n}\left(  r,X^{u}\left( r \right)  \right)Gu\left(  r\right)  \right] dr \to\int_{t}^{s}\left[  \psi\left( \nabla v\left(r,X^{u}\left(  r\right) \right)   G \right)  -\nabla v\left(  r,X^{u}\left( r \right)  \right)Gu\left(  r\right)  \right] dr 
\end{align*}
a.s and uniformly with respect to $s$.
So the sequence of local martingales $\zeta_{n}\left( s,X^{u}\left(
s\right)  \right)  $ converges a.s. uniformly with respect to $s,$ and
so also uniformly in probability, to the process
\begin{align*}
\zeta\left(  s,X^{u}\left( s\right)  \right)   &  :=v\left(s,X^{u}\left(  s\right)  \right)  -v(t,x)\\
&  +\int_{t}^{s}\left[  \psi\left( \nabla v\left(
r,X^{u}\left(  r\right)\right)  G \right)  -\nabla v\left(  r,X^{u}\left( r \right)  \right)Gu\left(  r\right)  \right]  dr.
\end{align*}
Consequently $\zeta\left(  s,X^{u}\left(  s\right)  \right)  $ is a local martingale, see \cite{K}. Moreover, by the representation (\ref{vsemimartingala}), it follows that%
\[
\zeta\left( s,X^{u}\left( s\right)  \right)  =\int_{t}^{s}Z\left(r\right)  dW\left( r\right)  +\int_{t}^{s}\left(  \nabla v\left(r,X^{u}\left(  r\right)  \right) G -Z\left(  r\right)\right)  u\left(  r\right)  dr.
\]
So for every $0\leq s\leq T$,
\[
\int_{t}^{s}\left(  \nabla v\left( r,X^{u}\left( r\right)  \right)G-Z\left( r\right)  \right)u\left(  r\right) dr=0,
\]
and this gives the identification of $Z\left(  s\right) $ with $\nabla v\left(  s,X^{u}(s)\right)G  $.
\newline With this identification in hands and taking into account \eqref{vsemimartingala}, we can rewrite the cost $J$ as 
\begin{equation*}
J(  t,x,u)  =v(t,x)+\E\int_{t}^{T}\left[ -\psi\left(\nabla v\left(  s,  X^{u}(s)  \right)G\right)    +\nabla v\left(  s,X^{u}(  s)  \right)G   u(  s)  +g(u(s))  \right] ds,
\end{equation*}
which gives the so called fundamental relation: by the definition of the hamiltonian function $\psi$ relation (\ref{fund_rel}) is proved. Moreover $J(t,x,u)=v(t,x)$ if and only of  f $u(s)=\Upsilon_{0}\left(  \nabla v(s X^{u}(s)G\right) $, $\mathbb{P}$--a.s. for almost every $0\leq t\leq s\leq T ,$, where $\Upsilon_0$ is measurable selection of $\Upsilon$ defined in \eqref{Gamma}.
\newline It is immediate to see that by the Girsanov Theorem equation (\ref{closed loop eq}) admits a  weak solution in mild sense, which turns out to be unique in law.
\qed


\begin{thebibliography}{99}

\bibitem{CM}A. Chojnowska-Michalik , \textit{Representation theorem for general stochastic delay equations.} Bull. Acad. Polon. Sci. Sér. Sci. Math. Astronom. Phys. 26 (1978), no. 7, 635--642. 

\bibitem{ContFou} R. Cont and D.-A. Fournié,  \textit{Functional Itô calculus and stochastic integral representation of martingales}. Ann. Probab. 41 (2013), no. 1, 109--133. 

\bibitem{CGRST} A. Cosso, S. Federico, F. Gozzi, M. Rosestolato, . Touzi, \textit{Path-dependent equations and viscosity solutions in infinite dimension}.
Ann. Probab. 46 (2018), no. 1, 126--174. 

\bibitem {DP1}G. Da Prato and J. Zabczyk,\textit{\ Stochastic equations in
infinite dimensions. Second Edition.}
Encyclopedia of Mathematics and its Applications 152,
Cambridge University Press, \textit{2014.}


\bibitem {DP2}G. Da Prato and J. Zabczyk,\textit{\ Ergodicity for
infinite-dimensional systems. }London Mathematical Society Note Series, 229,
Cambridge University Press, Cambridge, 1996.

\bibitem {DP3}G. Da Prato and J. Zabczyk,\textit{\ Second order partial
differential equations in Hilbert spaces}. London Mathematical Society Note
Series, 293, Cambridge University Press, Cambridge, 2002. 

\bibitem{DM}
{ M. C. Delfour and S. K. Mitter}, {\it Hereditary differential
systems with constant delay}, J. Differential Equations, 12
(1974), pp.~213--235.

\bibitem{Dup} B. Dupire,  \textit{Functional Itô calculus.} Quant. Finance 19 (2019), no. 5, 721--729. 

\bibitem{EKTZ}  I. Ekren, C. Keller, N. Touzi, Nizar, J. Zhang,  \textit{On viscosity solutions of path dependent PDEs.} Ann. Probab. 42 (2014), no. 1, 204--236.

\bibitem{FZ} F.Flandoli, G. Zanco, \textit{An infinite-dimensional approach to path-dependent Kolmogorov equations}. Ann. Probab. 44 (2016), no. 4, 2643--2693.

\bibitem{FRZ} F. Flandoli, F. Russo, G. Zanco,  \textit{Infinite-dimensional calculus under weak spatial regularity of the processes}. J. Theoret. Probab. 31 (2018), no. 2, 789–826. 

\bibitem{F} Fuhrman, M. \textit{Smoothing properties of nonlinear stochastic equations in Hilbert spaces}. NoDEA Nonlinear Differential Equations Appl. 3 (1996), no. 4, 445--464.

\bibitem {FT2}M. Fuhrman and G. Tessitore, \textit{Non linear Kolmogorov equations in infinite dimensional spaces: the backward stochastic differential equations approach and applications to optimal control. }Ann. Probab. 30 (2002), no. 3, 1397--1465.

\bibitem {FTGgrad}M. Fuhrman and G. Tessitore,\textit{Generalized directional gradients, backward stochastic differential equations and mild solutions of semilinear parabolic equations}.
Appl. Math. Optim. 51 (2005), no. 3, 279--332.

\bibitem{FMT} M. Fuhrman, F.  Masiero and G. Tessitore, \textit{Stochastic equations with delay: optimal control via BSDEs and regular solutions of Hamilton-Jacobi-Bellman equations}. SIAM J. Control Optim. 48 (2010), no. 7, 4624-4651. 

\bibitem {G1}F. Gozzi, \textit{Regularity of solutions of second order Hamilton-Jacobi equations in Hilbert spaces and applications to a control problem, }(1995) Comm Partial \ Differential Equations 20, pp. 775-826.

\bibitem{GoMa1}
{ F. Gozzi, F. Masiero}, {\it Stochastic optimal control with delay in the control I: Solving the HJB equation through partial smoothing.} SIAM J. Control Optim. 55 (2017), no. 5, pp. 2981-3012.


\bibitem{K}H. Kunita, \textit{Stochastic differential equations and stochastic
flows of diffeomorphisms.} \'{E}cole d'\'{e}t\'{e} de probabilit\'{e}s de
Saint-Flour, XII-1982, 143--303, Lecture Notes in Math., 1097, Springer,
Berlin, 1984.


\bibitem{Mas0} F. Masiero, \textit{Semilinear Kolmogorov equations and applications to stochastic optimal control,} PhD Thesis, 2004.

\bibitem {Mas}F. Masiero, \textit{Semilinear Kolmogorov equations and applications to stochastic optimal control}.  Appl. Math. Optim.  51  (2005),  no. 2, 201-250.

 \bibitem{Mas1} F. Masiero, \textit{Regularizing properties for transition semigroups and semilinear parabolic equations in Banach spaces.} Electron. J. Probab.  12  (2007), no. 13, 387--419.
 
 \bibitem{MOTZ} F. Masiero, C. Orrieri, G. Tessitore, G. Zanco, \textit{Semilinear {K}olmogorov equations on the space of continuous functions via {BSDE}s}, Stochastic Process. Appl.,136, (2021), 1--56.
	

 \bibitem{moha} { S.-E. A. Mohammed}, {\it Stochastic differential systems with memory: theory, examples and applications.}in: Stochastic Analysis and Related Topics, VI (Geilo, 1996), in: Progr. Probab., vol. 42, Birkh\"{a}user Boston, Boston,
MA, 1998, pp. 1-77.

\bibitem{RR} Z. Ren, M. Rosestolato, \textit{Viscosity solutions of path-dependent PDEs with randomized time}. SIAM J. Math. Anal. 52 (2020), no. 2, 1943--1979

\bibitem{RS}M. Rosestolato, A. Swiech, \textit{Partial regularity of viscosity solutions for a class of Kolmogorov equations arising from mathematical finance.} J. Differential Equations 262 (2017), no. 3, 1897--1930. 

\bibitem{W}
{G. F. Webb}, {\it Functional differential equations and
nonlinear semigroups in $L^{p}$-spaces}, J. Differential
Equations, 20 (1976), pp.~71--89.

\end{thebibliography}
\end{document}